\documentclass{birkjour}
\usepackage{amssymb}
\usepackage{orcidlink}
\newtheorem{thm}{Theorem}[section]
\newtheorem{cor}[thm]{Corollary}
\newtheorem{lem}[thm]{Lemma}
\newtheorem{prop}[thm]{Proposition}
\newtheorem{defn}[thm]{Definition}
\newtheorem{rem}[thm]{Remark}
\newtheorem*{ex}{Example}

\numberwithin{equation}{section}

\usepackage{refcount}

\hypersetup{
colorlinks=true,
linkcolor=blue,
citecolor=blue,
urlcolor=blue
}

\makeatletter
\renewcommand{\ref}[1]{\hyperref[#1]{\getrefnumber{#1}}}
\makeatother

\usepackage{tikz-cd}
\usetikzlibrary{backgrounds,fit,matrix,positioning,calc,through,chains}
\usetikzlibrary{arrows,shapes,decorations.pathmorphing,automata,petri}
\usepackage[boxsize=1.25em,centerboxes]{ytableau}

\makeatletter
\renewcommand{\thesection}{\arabic{section}}

\renewcommand{\theequation}{\thesection.\arabic{equation}}
\makeatother

\newcommand{\C}{{\mathbb C}}
\newcommand{\Z}{{\mathbb Z}}
\newcommand{\Q}{{\mathbb Q}}
\newcommand{\R}{{\mathbb R}}

\newcommand{\HH}{{\mathbb H}}

\newcommand{\OO}{{\mathbb O}}

\newcommand{\Lie}[1]{\mathfrak{#1}}

\newcommand{\mc}[1]{\mathcal{#1}}

\usepackage{xcolor}

\DeclareMathOperator{\End}{End}
\DeclareMathOperator{\GL}{GL}
\DeclareMathOperator{\Sp}{\mathbf{Sp}}

\begin{document}

%
%
%
%
%
%
%
%
%

\title[\small Cohomology of Homogeneous Spaces Associated with Composition Algebras]
 {(Equivariant) Cohomology of Homogeneous Spaces Associated with Composition Algebras}

\author[Mahir Bilen Can]{Mahir Bilen Can \orcidlink{orcid.org/0000-0002-0175-4897}}

\address{%
Department of Mathematics\\Tulane University\\
6823 St. Charles Avenue\\
New Orleans 70118--5698 LA\\USA}

\email{mahirbilencan@gmail.com}

\thanks{This work was completed with the support of our
\TeX-pert.}
\author{Ana Casimiro \orcidlink{orcid.org/0000-0001-6001-6803}}
\address{Department of Mathematics and Center for Mathematics and Applications\\ NOVA School of Sciences and Technology\\ Caparica 2829--516 \\ Portugal}
\email{amc@fct.unl.pt}

\author{Ferruh \"Ozbudak \orcidlink{orcid.org/0000-0002-1694-9283}}
\address{Faculty of Engineering and Natural Sciences \\ Sabanci University \\ Tuzla, Istanbul 34956 \\ Turkey}
\email{ferruh.ozbudak@sabanciuniv.edu}



\begin{abstract}
We study homogeneous spaces arising from embeddings of associative composition algebras of dimensions $2$ and $4$, with particular emphasis on the Hamiltonian quaternion algebra. For inclusions $D\hookrightarrow C$, we analyze the (equivariant) geometry and topology of the quotients $\mathbf{GL}(n,C)/\mathbf{GL}(n,D)$. In the quaternionic case, we prove that these spaces are equivariantly diffeomorphic to homogeneous vector bundles over compact symmetric spaces, and hence admit equivariant deformation retractions onto compact homogeneous models. This description determines their homotopy type and reduces the computation of their cohomological invariants to the compact setting.

Using this reduction together with Weyl group invariant theory, we compute the rational equivariant cohomology rings with respect to maximal tori and obtain explicit formulas for the Poincaré polynomials of the quaternionic homogeneous spaces. We further extend the analysis to noncompact symmetric spaces associated with Clifford groups. In this setting, we establish equivariant deformation retractions onto compact homogeneous spaces and derive explicit expressions for their Betti numbers.

These results provide a unified approach to the topology and equivariant cohomology of homogeneous spaces associated with composition and Clifford algebra embeddings.
\end{abstract}

\maketitle

\section{Introduction}

The starting point of our article is the following observation: Let $C$ be an associative composition algebra over a field $k$ of characteristic $\neq 2$, and let $D \subseteq C$ be a composition subalgebra. 
Then $C$ is naturally a finite-dimensional left $D$-vector space. 
The left regular action
\[
\lambda : C \to \End_D(C), \qquad \lambda(c)(x)=cx,
\]
is $D$-linear and injective. 
Consequently, $C$ embeds as a $k$-subalgebra of $\End_D(C)$ and therefore admits a faithful $D$-linear representation.
By a composition algebra over a field $k$, we mean a $k$-algebra $C$ equipped with a nondegenerate quadratic form $N: C \to k$ such that $N(xy) = N(x)N(y)$ for all $x, y \in C$. 
Examples of associative composition algebras include all fields, all quadratic algebras over fields, and all quaternion algebras over fields. 
The central question of our study is the following: What is the relationship between the groups of linear automorphisms of $C$ and $D$? More broadly, how can one describe the algebra (equivariant $K$-theory) and the geometry
(equivariant cohomology) of the homogeneous space measuring the difference between
$C$-linear and $D$-linear symmetries?
More precisely, for an associative composition algebra $C$ and a composition subalgebra
$D\subseteq C$, we consider the free left $C$-module $C^n$ together with its natural
left $D$-submodule $D^n\subset C^n$. The corresponding homogeneous space is
\[
\mathbf{GL}(n,C)\big/\mathbf{GL}(n,D),
\]
which measures the difference between $C$-linear and $D$-linear automorphisms 
of the free module $C^n$.
In the case where $C$ is obtained from $D$ by the Cayley--Dickson doubling 
construction, this homogeneous space is symmetric (see Proposition~\ref{intro:P1}). 
In general, however, the space need not be symmetric.
Our goal is to investigate the geometric invariants of these 
homogeneous spaces, in particular their cohomology and  equivariant cohomology.

To begin, we present our cohomological results with a specific focus on associative composition algebras defined over the real numbers. 
These algebras, in dimensions 2 and 4, up to isomorphism, are given as in Table~\ref{T:list1}. 
\begin{table}[h]
	\centering
	\begin{tabular}{|c|c|}
		\hline
		\textbf{Real Non-split Composition Algebras} & \textbf{Real Split Composition Algebras} \\
		\hline
		$\C$ & $\C_s := \R \oplus \R$ \\
		\hline
		$\HH = $ Hamilton's quaternions & $\HH_s := $ real $2\times 2$ matrices \\
		\hline
	\end{tabular}
	\caption{Associative real composition algebras in dimensions 2 and 4.}
	\label{T:list1}
\end{table}

Each algebra in the first row of Table~\ref{T:list1} has a natural embedding into the one directly below it, making the latter an algebra over the former. 
Additionally, there is an $\R$-algebra map from $\C$ into $\HH_s$.
For $C \in \{ \C,\C_s,\HH,\HH_s\}$, the $n$-fold product $C^n$ is a free $C$-module with respect to the right $C$-action. 
Of course, $\mathbf{GL}(n,C)$ is a real algebraic group for $C\in \{\C,\HH\}$.
We observe that $\mathbf{GL}(n,C)$, is a real algebraic group for $C\in \{\C_s,\HH_s\}$ as well.
Let $D$ be another associative composition algebra from $\{ \C,\C_s,\HH,\HH_s\}$.
Notice that if there is an injective $\R$-algebra homomorphism $D\hookrightarrow C$, then there is a corresponding injective homomorphism 
of real algebraic groups, $\mathbf{GL}(n,D)\hookrightarrow \mathbf{GL}(n,C)$.

For cohomological computations we pass to maximal compact subgroups. 

Let $\mathbf{U}(1)$ denote the unitary group of complex numbers of unit length.
Let $S_n$ denote the symmetric group on $n$ letters.
Let $T=\mathbf{U}(1)^n$. Then
\[
\mathrm{Lie}(T)\cong \R^n,
\]
and the symmetric algebra on its dual is
\[
S(\Lie(T)^*) \cong \R[x_1,\dots,x_n].
\]
We write
\[
S:=\Q[x_1,\dots,x_n],\qquad \deg x_i=2,
\]
so that $S\cong H_T^*(pt;\Q)$.
The Weyl group of $\mathbf{SL}(n,\R)$ with respect to $T$ is
\[
W(\mathbf{SL}(n,\R))\cong S_n.
\]
We denote by $W_n$ the Weyl group of $\mathbf{SO}(n,\R)$ (of type $B_m$ or $D_m$ according to the parity of $n$). For a compact Lie group $G$ containing $T$ as a maximal torus,
denote by $W(G)$ its Weyl group with respect to $T$.

In particular,
\[
W(\mathbf{U}(n)) \cong S_n,
\]
while $W(\mathbf{Sp}(n))$ and $W(\mathbf{SO}(2n))$ are of type $C_n$ and 
$D_n$, respectively.

In this notation, the first main result of our article is the following theorem.

\begin{thm}\label{intro:1}
Let $D$ and $C$ be composition algebras in $\{\C,\C_s,\HH,\HH_s\}$
such that there exists an injective algebra homomorphism $D\hookrightarrow C$.

Then the following assertions hold.

\begin{enumerate}

\item[(1)] If $(C,D)=(\HH,\C)$, then the Poincar\'e polynomial of the homogeneous space $\mc{Y}_n:= \mathbf{GL}(n,\HH)/\mathbf{GL}(n,\C)$ is
$
P(\mc{Y}_n,t)=\prod_{i=1}^{n}(1+t^{2i}).
$
Moreover, there is an isomorphism of graded $S$-algebras
\[
H_T^*(\mc{Y}_n;\Q)
\;\cong\;
S \otimes_{S^{W(\mathbf{Sp}(n))}}
S^{W(\mathbf{U}(1)\cdot \mathbf{SU}(n))}.
\]

\item[(2)] If $(C,D)=(\HH_s,\C)$, then the Poincar\'e polynomial of the homogeneous space $\mc{Z}_n:= \mathbf{GL}(n,\HH_s)/\mathbf{GL}(n,\C)$ is 
$
P(\mc{Z}_n,t)=\prod_{i=1}^{n-1}(1+t^{2i}).
$
Moreover, there is an isomorphism of graded $S$-algebras
\[
H_T^*(\mc{Z}_n;\Q)
\;\cong\;
S \otimes_{S^{W(\mathbf{SO}(2n))}}
S^{W(\mathbf{U}(1)\cdot \mathbf{SU}(n))}.
\]

\item[(3)] If $(C,D)=(\HH_s,\C_s)$, then the Poincar\'e polynomial of the corresponding homogeneous space $\mc{W}_n:= \mathbf{GL}(n,\HH_s)/\mathbf{GL}(n,\C_s)$ is
$
P(\mc{W}_n,t)=
\begin{bmatrix}
n \\ \lfloor n/2\rfloor
\end{bmatrix}_{t^2}
$ (Gaussian binomial coefficient).
Moreover, there is an isomorphism of graded $S$-algebras
\[
H_T^*(\mc{W}_n;\Q)
\;\cong\;
S \otimes_{S^{W(\mathbf{SO}(2n))}}
S^{W(\mathbf{SO}(n)\times \mathbf{SO}(n))}.
\]

\end{enumerate}
All cohomology rings are taken with rational coefficients.
\end{thm}
\medskip

Our previous theorem serves as a warm-up for a result that provides much more algebraic information on a general family of homogeneous spaces. 
Let $G$ be a connected reductive algebraic group defined over a field $k$. 
Let $\sigma : G\to G$ be an automorphism of $G$ such that $\sigma^2 = id$.
Let $G^\sigma  :=\{g\in G \mid \sigma(g) = g\}$.
If a closed subgroup $H \subseteq G$ satisfies $G^\sigma \subseteq H \subseteq N_G(G^\sigma)$, where $N_G(G^\sigma)$ is the normalizer of $G^\sigma$ in $G$,
then the homogeneous space $G/H$ is called a {\em symmetric space}. 
In this terminology, a simple but key proposition that we based our work is the following observation.

\begin{prop}  \label{intro:P1}
	Let $C$ be an associative composition algebra over a field $k$. 
	Let $D$ be a composition subalgebra of $C$ such that $C$ is obtained from $D$ by the Cayley-Dickson doubling process.
	For $n\geq 1$, let $G$ (resp. $H$) denote $\mathbf{GL}(n,C)$ (resp. $\mathbf{GL}(n,D)$).
	Then, after base change to a splitting field $S$ of $C$,
the homogeneous space
\[
G_S/H_S
\]
is a symmetric space.
\end{prop}

After base change to a splitting field $S$ of $C$, 
the algebra $C_S$ becomes isomorphic to a matrix algebra, and the 
conjugation involution arising from the Cayley–Dickson construction 
induces an algebraic involution of $G_S$ whose fixed-point subgroup 
contains $H_S$. This explains the symmetric structure.

For background on composition algebras and their construction using the Cayley-Dickson doubling process, we recommend the monograph~\cite{VeldkampSpringer}.

We are now ready to give a brief overview of our paper and explain its structure.
In Section~\ref{S:Preliminaries}, we discuss our motivational composition algebras in more detail. 
In particular, we introduce quaternions since they are our four dimensional composition algebras. 
In Section~\ref{S:Endomorphisms}, the endomorphism ring of the free $C$-modules (of finite rank) are introduced. 
Therein, we not only review some well-known results from the literature regarding the defining representations of the general linear groups $\mathbf{GL}(n,C)$, where $C$ is an associative composition algebra, but also introduce our symmetric spaces and prove our 
Proposition~\ref{intro:P1}, showing that the Cayley–Dickson doubling construction naturally 
leads to symmetric spaces after base change.
In Section~\ref{S:Quaternionic}, we focus on the geometry of our symmetric spaces. 
Among other things, we prove our theorem on the $\mathbf{U}(1)$-equivariant cohomology rings, that is, Theorem~\ref{intro:1}. 
We close our paper in Section~\ref{S:Final} after explaining an important special case of our opening question in relation with the Clifford algebras and spin groups.

\section{Preliminaries}\label{S:Preliminaries}

It is a well-known theorem of Hurwitz~\cite{Hurwitz1922} that every normed associative division algebra over the reals is isomorphic to either the real numbers, the algebra of complex numbers, or Hamilton's quaternions, denoted $\HH$.
Hurwitz showed that, dropping the associativity requirement adds only one more algebra to this list, namely the Cayley algebra, denoted by $\OO$. Although every nonzero element of $\OO$ has a multiplicative inverse, due to nonassociativity, $\OO$ is not a (skew-)field.  
After Jacobson's work~\cite{Jacobson1958}, we know that all of these normed division algebras are special examples of the composition algebras that we introduced at the beginning of this article. 
A composition algebra $C$ over $k$ is said to be {\em split} if it contains zero divisors. 
In each of the dimensions 2, 4, and 8, there is, up to isomorphism, exactly 
one split composition algebra, and these are the only composition algebras containing zero divisors. 
In fact, according to Jacobson, every two dimensional split composition algebra is isomorphic to $k \oplus k$, every four dimensional split composition algebra is isomorphic to the algebra of $2\times 2$ matrices with entries from $k$,
and every eight dimensional split composition algebra is isomorphic to the composition algebra obtained from $2\times 2$ matrices by applying the Cayley-Dickson doubling process. 
This last split composition algebra over the reals is called the {\em split Cayley} algebra, which we denote by $\OO_s$. 
Now we have a complete, two column list of all real composition algebras, up to isomorphism, in dimensions 2, 4, and 8 over reals:
\begin{table}[h]
	\centering
	\begin{tabular}{|c|c|}
		\hline
		\textbf{Normed Division Algebras} & \textbf{Split Composition Algebras} \\
		\hline
		$\C$ & $\C_s \cong \R \oplus \R$ \\
		\hline
		$\HH$ & $\HH_s \cong \mathbf{Mat}(2,\R)$ \\
		\hline
		$\OO$ & $\OO_s$ \\
		\hline
	\end{tabular}
	\caption{Real composition algebras in dimensions 2,4, and 8.}
	\label{T:list}
\end{table}

\begin{rem}
	The goal of our paper is to study the geometry of the homogeneous spaces associated with the entries of Table~\ref{T:list1}. 
	In particular, we will consider the general linear groups defined over associative composition algebras of dimensions 2 and 4.
	The reason for why we do not handle linear groups over octonions (Cayley algebras, split or not) is because they are nonassociative;
	the study of automorphism groups of their free modules requires a different line of investigation that does not align with the main geometric results of the present paper. 
	For this reason, hereafter we will not discuss Cayley algebras.
\end{rem}

\subsection{Inclusions of algebras.}

An associative composition algebra $C$ over $k$, a field of characteristic $\neq 2$, such that $\dim_k C = 4$ is a quaternion algebra. 
Quaternion algebras are examples of central simple $k$-algebras. 
The purpose of this subsection is to review the basics of such composition algebras. 
An informative and accessible account of this theory (quaternion algebras) can be found in the unpublished notes by Keith Conrad~\cite{KConrad}. 
The proofs of our assertions in this subsection are either contained in~\cite{KConrad} or can be easily derived from the information therein.
\medskip

A {\em quaternion algebra over $k$} is a four dimensional $k$-algebra $C$ having a vector space basis $\{1,u,v,w\} \subseteq C$ 
such that 
\begin{enumerate}
	\item $\{a,b\} \subset k^*$, where $a=u^2$ and $b=v^2$, 
	\item $w= uv = -vu$, 
\end{enumerate}
In this case, we write $(a,b)_k$ for $C$. 
Note that there is a natural isomorphism, 
$$
(a,b)_k\cong (b,a)_k.
$$

Let $z\in (a,b)_k$ be given by $z=  x_0 +x_1u + x_2 v + x_3 w$ where $x_i\in k$ ($i\in \{0,1,2,3\}$). 
If $x_0 = 0$, then $z$ is said to be a {\em pure quaternion}. 
The {\em conjugate} and the {\em norm} of $z$ are defined by 
$$
\overline{z} := x_0 - x_1u - x_2 v - x_3 w\quad\text{and}\quad N(z) := z \overline{z}.
$$
A quaternion $z\in (a,b)_k$ has a 2-sided inverse if and only if $N(z)\neq 0$. 
In this case, the inverse of $z$ is given by $\overline{z}/N(z)$. 
Assuming that ${\rm char}\ k\neq 2$, the center of a quaternion algebra $(a,b)_k$ is given by $k$. 
If ${\rm char}\ k= 2$, then $(a,b)_k$ is commutative. 
In the rest of this article, we assume that the base field of any of our composition algebras is not two, ${\rm char}\ k\neq 2$.
\medskip

Recall that a composition algebra $C$ over $k$ is called split 
if it contains zero divisors. 
In dimension four, this condition is equivalent to 
$C \cong \mathbf{Mat}(2,k)$; thus a quaternion algebra is split 
if and only if it is isomorphic to the matrix algebra $\mathbf{Mat}(2,k)$. 
For a field extension $k\subseteq F$, a quaternion algebra $C$ over $k$ is said to {\em split over $F$} if after base change to $F$ it becomes a split quaternion algebra over $k$. 
In other words, $C$ splits over $F$ if $C\otimes_k F \cong \mathbf{Mat}(2,F)$.
Note that, over an algebraically closed field every quaternion algebra is split. 
Equivalently, if $\overline{k}$ denotes an algebraic closure of $k$ and $C_{\overline{k}}$ denotes the quaternion algebra $C\otimes_k \overline{k}$, then we have 
\begin{align}\label{A:basechangeforC}
	C_{\overline{k}} \cong \mathbf{Mat}(2,\overline{k}).
\end{align}
\medskip

Let $C$ be a quaternion algebra over $k$. 
Since $C$ is four-dimensional over $k$, 
the left regular representation embeds $C$ into $\mathbf{Mat}(4,k)$.
If $C$ is split, then $C \cong \mathbf{Mat}(2,k)$, 
and hence it admits a faithful $2 \times 2$ matrix representation over $k$. 
Let us assume that $C$ is a non-split quaternion algebra. 
We proceed to construct a faithful matrix representation of $C$. 
Let $L\subset (a,b)_k$ be a maximal subfield of $C$. 
Then it is necessarily true that $k\subsetneq L \subsetneq C$ and that $L$ is a two dimensional sub-composition algebra of $C$.
In fact, $L$ is a quadratic field extension of $k$: it is given by $L= k(\sqrt{a})$ for some $a\in k^*$.
Then we have an involutive field automorphism, 
\begin{align*}
	\tau : L &\longrightarrow L \\
	x+y \sqrt{a} &\longmapsto x-y\sqrt{a}.
\end{align*}
Since $C$ is a two dimensional $L$-algebra, we have a decomposition 
\begin{align*}
	C = 1 L \oplus v L,
\end{align*}
for some $v\in C^*$ such that $z v = v \tau(z)$ for every $z\in L$. 
Let $u$ denote $\sqrt{a}$. 
Then we have $uv = - vu$. 
Notice also that $zv^2 = zvv = v \tau(z) v = v^2 \tau^2(z)=v^2z$.
It follows that $v^2$ commutes with every element of $C$.
Since $C$ is a central $k$-algebra, this means that $v^2 \in k$.
Let $b:=v^2$. 
It is now evident that $C=(a,b)_k$. 
\medskip

We consider the assignments 
$$
1 \mapsto \begin{bmatrix} 1 & 0 \\ 0 & 1 \end{bmatrix}, \quad u \mapsto \begin{bmatrix} \sqrt{a} & 0 \\ 0 & -\sqrt{a} \end{bmatrix}, \quad 
v \mapsto \begin{bmatrix} 0& -1 \\ -b & 0 \end{bmatrix}, \quad w \mapsto \begin{bmatrix} 0 & -\sqrt{a} \\ \sqrt{a}b & 0 \end{bmatrix}.
$$
Extended to $(a,b)_k$ by $k$-linearity, these assignments define a faithful linear representation of $(a,b)_k$ that we sought after. 
We call it the {\em symplectic representation of $(a,b)_k$}.
We denote it by 
\begin{align}\label{A:symplecticrep}
	\sigma : (a,b)_k \longrightarrow \mathbf{Mat}(2,L).
\end{align}
Then the matrix form of the ``Cayley-Dickson doubling process'' takes the following shape:  
\begin{align}\label{A:CDdoubling1}
	\sigma((a,b)_k) & =\left\{ \begin{bmatrix}  x_1+ y_1\sqrt{a} & -(x_2- y_2\sqrt{a}) \\ -b(x_2+ y_2\sqrt{a}) & x_1-y_1\sqrt{a}  \end{bmatrix} \mid x_1,x_2,y_1,y_2\in k\right\}.
\end{align}
In particular, for $b= -1$, we see that the image of $\sigma$ is given by 
\begin{align}\label{A:imageofsigma}
	\sigma((a,-1)_k) = \left\{ \begin{bmatrix}  x_1+ y_1\sqrt{a} & -(x_2- y_2\sqrt{a}) \\ x_2+ y_2\sqrt{a} & x_1-y_1\sqrt{a}  \end{bmatrix} \mid x_1,x_2,y_1,y_2\in k\right\}.
\end{align}
Such a quaternion algebra is split if and only if $a$ is a sum of two squares in $k$.
A proof of this fact can be found in~\cite[Corollary 4.23]{KConrad}.

Finally, it is worth mentioning that the maximal subfield $L$ of $(a,b)_k$ is mapped under $\sigma$ to the diagonal subalgebra, 
\begin{align}\label{A:CDdoubling12}
	L\cong \left\{ \begin{bmatrix}  x_1+ y_1\sqrt{a} &0 \\0  & x_1-y_1\sqrt{a}  \end{bmatrix} \mid x_1,x_2,y_1,y_2\in k\right\},
\end{align}
where the $L$-algebra structure on $(a,b)_k$ is given by the right multiplication action of (\ref{A:CDdoubling12}) on (\ref{A:CDdoubling1}).

\begin{ex}
	The Hamilton quaternion algebra $\HH$ is given by $(-1,-1)_\R$.  
	In the literature, the notation for the basis $\{1,u,v,w\}$ is usually given by $\{1,i,j,l\}$.
	Hence, the elements of $\HH$ are given by $a_0  + a_1i+ a_2j + a_3 l$ where $a_0,a_1,a_2,a_3\in \R$ and $i,j,l$ satisfy the following relations:
	$$
	i^2 = j^2 = l^2 = -1 ,\quad ij =l = -ji, \quad jl = i = - lj,\quad li = j = - il.
	$$
	Notice that the real subalgebra generated by 1 and $i$ is isomorphic to the field of complex numbers, $\C$. 
	This copy of $\C$ in $\HH$ is a maximal subfield of $\HH$. 
	It plays the role of $L$ for a general non-split quaternion algebra. 
	We have the decomposition  
	\begin{align}\label{A:CD1}
		\HH =  1\C \oplus  j \C. 
	\end{align}
	The symplectic representation of $\HH$, as a special case of (\ref{A:imageofsigma}), is given by 
	\begin{align*}
		\sigma : \HH &\longrightarrow \mathbf{Mat}(2,\C) \\
		z+ wj &\longmapsto \begin{bmatrix} z & -\overline{w} \\ w & \overline{z} \end{bmatrix}.
	\end{align*}
	The Cayley-Dickson decomposition (\ref{A:CD1}) in terms of matrices is given by 
	\begin{align*}
		\underbrace{\begin{bmatrix} z & -\overline{w} \\ w & \overline{z} \end{bmatrix}}_{1 z+ j w \in \HH} = 
		\underbrace{\begin{bmatrix} 1& 0 \\ 0 & 1 \end{bmatrix}}_{1} \underbrace{\begin{bmatrix} z &0 \\ 0 & \overline{z} \end{bmatrix}}_{z}+
		\underbrace{\begin{bmatrix} 0 & -1\\ 1 & 0 \end{bmatrix}}_{j}\underbrace{\begin{bmatrix} w & 0\\ 0 &\overline{w} \end{bmatrix}}_{w}.
	\end{align*}
	
\end{ex}

\begin{ex}\label{E:splitcase}
	In this example we consider split quaternion algebras over a field $k$. 
	These are precisely the four dimensional split associative composition algebras over $k$. 
	In particular, they are all isomorphic to the matrix algebra $\mathbf{Mat}(2,k)$. 
	Let us explain the Cayley-Dickson doubling process for this matrix algebra.
	To this end, we will use the following basis (recall ${\rm char}\ k \neq 2$):
	\begin{align}\label{E:gens}
		\mathbf{1}_2:= \begin{bmatrix} 1 & 0 \\ 0 & 1 \end{bmatrix}, \quad u:= \begin{bmatrix} 1 & 0 \\ 0 & -1 \end{bmatrix}, \quad 
		v := \begin{bmatrix} 0& -1 \\ 1 & 0 \end{bmatrix}, \quad w := \begin{bmatrix} 0 & -1 \\ -1 & 0 \end{bmatrix}.
	\end{align}
	By using this basis, it is easy to see that $\mathbf{Mat}(2,k)$ is isomorphic to $(1,-1)_k$. 
	Let $D$ denote the diagonal subalgebra of $\mathbf{Mat}(2,k)$ generated by $\mathbf{1}_2$ and $u$.
	Then $D$ is not a maximal subfield but a two dimensional sub-composition algebra of $C$.
	The Cayley-Dickson doubling process in this case is then given by 
	\begin{align}\label{A:CDdoublining2}
		\mathbf{Mat}(2,k) = \mathbf{1}_2 D \oplus w D.
	\end{align}
\end{ex}

\bigskip

We continue Example~\ref{E:splitcase} by constructing a maximal subfield.
Let $L$ denote the subset 
\begin{align}\label{A:CintoM2}
	L:= \left\{ \begin{bmatrix} x & -y \\ y & x \end{bmatrix} \mid x,y\in k \right\}.
\end{align}
Evidently, $L$ is closed under matrix addition. 
Furthermore, every nonzero element of $L$ has an inverse in $L$.
Indeed, for every $(x,y)\in k\times k$ such that $(x,y)\neq (0,0)$, we have 
$$
\begin{bmatrix} x & -y \\ y & x \end{bmatrix} ^{-1} = \frac{1}{x^2+y^2}\begin{bmatrix} x & y \\ -y & x \end{bmatrix} \in L.
$$
It is easy to check that $L$ is commutative:
$$
\begin{bmatrix} x & -y \\ y & x \end{bmatrix} \begin{bmatrix} z & -w \\ w & z \end{bmatrix} 
= \begin{bmatrix} xz -yw & -xw-yz \\ yz+xw & -yw +xz \end{bmatrix}
= \begin{bmatrix} z & -w \\ w & z \end{bmatrix} \begin{bmatrix} x & -y \\ y & x \end{bmatrix}.
$$
The remaining field properties are easy to check as well. 
Hence, $L$ is a subfield of $\mathbf{Mat}(2,k)$.
In fact, if $\sqrt{-1}\notin k$, then we can identify $L$ with the quadratic field extension $k[t]/(t^2+1)$, where $t$ corresponds to the matrix 
$\begin{bmatrix} 0 & -1 \\ 1 & 0 \end{bmatrix}$ and 1 corresponds to the identity matrix $\begin{bmatrix} 1 & 0 \\ 0 & 1 \end{bmatrix}$. 
Hence, the element $x+yt \in k[t]/(t^2+1)$ is represented by the matrix $\begin{bmatrix} x & -y \\ y & x \end{bmatrix}$. 
\medskip

We close our preparatory section by a remark that indicates when we can expect to find a homogeneous space from a pair of composition algebras.

\begin{rem}
	If $D$ is a split composition algebra, then it is not possible to embed it into any nonsplit composition algebra $C$ since the former algebra contains zero divisors. 
\end{rem}

\section{Endomorphisms and the Symmetric Space of a Maximal Subalgebra}\label{S:Endomorphisms}

Let $C$ be an associative composition algebra over a field $k$. 
In this section, following the notation of the previous subsection, we discuss the ring of endomorphisms of the free $C$-module of rank $n$.

We begin with the nonsplit case. 
Let $C:=(a,b)_k$. 
Let $\{1,u,v,uv\}$ be the basis such that $u^2=a$ and $v^2 = b$.
We use the general form of our conjugation to define the ``symplectic representation'' of the ring $\mathbf{Mat}(n,C)$. Let $L = k(\sqrt{a})$ be a maximal subfield of $C$, 
and let $\tau$ denote the nontrivial $k$-automorphism of $L$.
To this end, for $Z = (z_{ij})_{i,j=1}^n \in \mathbf{Mat}(n,C)$, we build two $n\times n$ matrices $A$ and $B$ with entries from $L$ as follows. 
Let $\sigma$ denote the symplectic representation as in (\ref{A:symplecticrep}). 
Let $i,j\in \{1,\dots, n\}$. 
We write the $(i,j)$-th entry of $Z$ as follows 
$$
z_{ij} = 1 x_{ij} + v y_{ij}, \quad \text{ where $x_{ij},y_{ij} \in L$. }
$$
Then, we have the decomposition 
$$
Z = X + v  Y,
$$
where $X=(x_{ij})_{i,j=1}^n\ ,Y=(y_{ij})_{i,j=1}^n \in \mathbf{Mat}(n,L)$. 
We set 
$$
\tau(Y):=(\tau(y_{ij}))_{i,j=1}^n.
$$ 
In this notation, we have the following generalization of (\ref{A:imageofsigma}):
\begin{align}\label{A:definingrepgeneral}
	\sigma: \mathbf{Mat}(n,C) &\longrightarrow \mathbf{Mat}(2n, L) \notag \\
	X+ vY &\longmapsto \begin{bmatrix} X & -\tau(Y) \\ -bY & \tau(X)\end{bmatrix}.
\end{align}
We refer to this representation as the {\em symplectic representation of $\mathbf{Mat}(n,C)$ over $L$}. 
We call the restriction of (\ref{A:definingrepgeneral}) to the group of invertible $C$-linear endomorphisms of $C^n$, that is, $\mathbf{GL}(n,C)$, the {\em symplectic representation of $\mathbf{GL}(n,C)$ over $L$}, as well.
We proceed to discuss $\mathbf{GL}(n,C)$ in more detail.

\medskip

\textbf{Over non-split quaternions.}

We begin with the main motivating example of our paper. 
  \begin{ex}\label{E:HamiltonionGL}
    Let $k=\R$ and $a=b=-1$. 
Then $C=(a,b)_k=\HH$. 
Every element $Z \in \mathbf{GL}(n,\HH)$ can be written uniquely in the form
\[
Z = A + jB, \qquad A,B \in \mathbf{Mat}(n,\C),
\]
via the identification $\HH=\C \oplus j\C$.
Then the symplectic representation of $\mathbf{GL}(n,\HH)$ over $\C$ is given by  
	\begin{align}\label{A:definingrep}
		\sigma: \mathbf{GL}(n,\HH) &\longrightarrow \mathbf{GL}(2n, \C) \notag \\
		A+ jB &\longrightarrow \begin{bmatrix} A & -\overline{B} \\ B & \overline{A}\end{bmatrix}.
	\end{align}
realizes $\mathbf{GL}(n,\HH)$ as a real algebraic subgroup of $\mathbf{GL}(2n,\C)$.
To compute its dimension, we argue directly. 
Since $\dim_\R \HH = 4$, the space
\[
\mathbf{Mat}(n,\HH)
\]
is an $\R$-vector space of dimension $4n^2$. 
Moreover,
\[
\mathbf{Gl}(n,\HH)=\mathbf{Aut}_{\HH}(\HH^n)
\]
is the open subset of invertible elements in $\mathbf{Mat}(n,\HH)$.
Therefore, as a real algebraic variety,
\[
\dim_\R \mathbf{Gl}(n,\HH)=\dim_\R \mathbf{Mat}(n,\HH)=4n^2.
\]
\end{ex}
\medskip

The analysis of the general linear groups with entries from a nonsplit quaternion algebra is not very different from the analysis for $\mathbf{GL}(n,\HH)$. 
Indeed, by using the symplectic representation of $\mathbf{GL}(n,C)$ over $L$ we see that 
\begin{align}\label{A:definingrepforGL(nC)}
	\mathbf{GL}(n,C) \cong \left\{  \begin{bmatrix} X & -\tau(Y) \\ -bY & \tau(X)\end{bmatrix} \in \mathbf{GL}(2n,L) \mid X,Y\in \mathbf{Mat}(n,L) \right\}.
\end{align}
Clearly, the dimension count argument for $\mathbf{GL}(n,\HH)$ extends to the case of $\mathbf{GL}(n,C)$.
As a $k$-variety the dimension of $\mathbf{GL}(n,C)$ is given by 
\begin{align}\label{A:nonsplitGLndim}
	\dim \mathbf{GL}(n,C) = 4n^2\qquad (\text{as a $k$-variety}).
\end{align} 
This follows directly from the fact that $\mathbf{GL}(n,C)$ is dense in $\mathbf{Mat}(n,C)$, the faithful representation of $\mathbf{Mat}(n,C)$ 
described in (\ref{A:definingrepgeneral}), and $\dim_k L = 2$. 
\medskip

\textbf{Over split quaternions.}

Next, we consider the ring $\mathbf{Mat}(n,C)$, where $C = \mathbf{Mat}(2,k)$. 
Then every element of $\mathbf{Mat}(n, \mathbf{Mat}(2,k))$ can be viewed as a $2n\times 2n$ matrix with entries from $k$. 
Let 
\[
Z = \begin{pmatrix}
	Z_{11} & Z_{12} & \cdots & Z_{1n} \\
	Z_{21} & Z_{22} & \cdots & Z_{2n} \\
	\vdots & \vdots & \ddots & \vdots \\
	Z_{n1} & Z_{n2} & \cdots & Z_{nn} \\
\end{pmatrix}
\quad\text{ and }\quad  
Y = \begin{pmatrix}
	Y_{11} & Y_{12} & \cdots & Y_{1n} \\
	Y_{21} & Y_{22} & \cdots & Y_{2n} \\
	\vdots & \vdots & \ddots & \vdots \\
	Y_{n1} & Y_{n2} & \cdots & Y_{nn} \\
\end{pmatrix}
\]
be two matrices from $\mathbf{Mat}(2n, k)$ written in block form, where each of the blocks is a $2\times 2$ matrix from $\mathbf{Mat}(2,k)$. 
Since the matrix product $ZY$ can be found by the block-matrix multiplication, we see that the set theoretic identification
$$
\mathbf{Mat}(n, \mathbf{Mat}(2,k)) \cong \mathbf{Mat}(2n,k)
$$
is in fact an $k$-algebra isomorphism. 
It follows that the invertible elements of $\mathbf{Mat}(n,C)$ is nothing but the general linear group of $2n\times 2n$ matrices with entries from $k$,
\begin{align}\label{A:importantstep}
	\mathbf{GL}(n, \mathbf{Mat}(2,k)) = \mathbf{GL}( 2n, k).
\end{align}
In particular, we have 
\begin{align}\label{A:splitGLndim}
	\dim \mathbf{GL}(n,C) = 4n^2\qquad (\text{as a $k$-variety}).
\end{align} 
We collect these observations in a proposition.

\begin{prop}\label{P:geometricpoints}
	Let $C$ be a quaternion algebra over a field $k$. 
	Let $n\geq 1$. 
	Then as a $k$-variety, the dimension of the general linear group of invertible $C$-linear transformations on $C^n$ is $4n^2$. 
	Furthermore, the group of geometric points of $\mathbf{GL}(n,C)$ is isomorphic to $\mathbf{GL}(2n,\overline{k})$, where 
	$\overline{k}$ is an algebraic closure of $k$. 
\end{prop}
\begin{proof}
	The split case is recorded in (\ref{A:splitGLndim}).
	The non-split case is recorded in (\ref{A:nonsplitGLndim}).
	These remarks take care of the first part of the proof of our proposition. 
	
	For our second assertion, we see from (\ref{A:importantstep}) that, if $C$ is split, then $\mathbf{GL}(n,C) \otimes_k \overline{k}$ 
	is isomorphic to $\mathbf{GL}(2n,\overline{k})$. 
	If $C$ is nonsplit, then the group of geometric points of $\mathbf{GL}(n,C)$ is given by 
	$$
	\mathbf{GL}(n,C)\otimes_k \overline{k} \cong \mathbf{GL}(n,C\otimes_k \overline{k}).
	$$
	But we know from (\ref{A:basechangeforC}) that $C\otimes_k \overline{k}$ is isomorphic to the split quaternion algebra $\mathbf{Mat}(2,\overline{k})$. 
	Hence, our problem reduces to the split case.
	This finishes the proof of our proposition. 
\end{proof}
\bigskip

In the rest of this subsection we describe the endomorphism ring of the two-dimensional composition algebras over a field. 
Let $D$ be a such a composition algebra over $k$. 
\medskip

First, we assume that $D$ is a nonsplit composition algebra over $k$ such that $\dim_k D =2$. 
Then $D$ is a quadratic field extension, $D= k(\sqrt{a})$ for some $a\in k$ such that $a^2\notin k$.
In this case, there is nothing mysterious. 
Indeed, the familiar group $\mathbf{GL}(n,D)$ is the group of $D$-linear automorphisms of the vector space $D^n$. 

There is a representation of $\mathbf{GL}(n,D)$ on $k^{2n}$. 
It is given by the restriction of the faithful ring representation 
\begin{align}\label{A:Cmonoidrep}
	\rho: \mathbf{Mat}(n,D) &\longrightarrow \mathbf{Mat}(2n, k) \notag \\
	A+ \sqrt{a} B &\longrightarrow \begin{bmatrix} A & aB \\ B & A\end{bmatrix}.
\end{align}
Hence, we have the identification 
\begin{align*}
	\mathbf{GL}(n,D) \cong \left\{ \begin{bmatrix} A & aB \\ B & A\end{bmatrix} \in \mathbf{GL}(2n,k) \mid A,B \in \mathbf{Mat}(n,k) \right\}.
\end{align*}
It follows that the dimension of $\mathbf{GL}(n,D)$ as a $k$-group is $2n^2$. 
\medskip

Finally, we consider the two dimensional split composition algebra, $D = k\oplus k$. 
The arguments that we used in the previous case are applicable to this case after we identify $D$ with the algebra of diagonal matrices in $\mathbf{Mat}(2,D)$. 
In particular, we have 
\small\[
\mathbf{GL}(n,D) = 
\left\{ \begin{pmatrix}
	Z_{11} & Z_{12} & \cdots & Z_{1n} \\
	Z_{21} & Z_{22} & \cdots & Z_{2n} \\
	\vdots & \vdots & \ddots & \vdots \\
	Z_{n1} & Z_{n2} & \cdots & Z_{nn} \\
\end{pmatrix}\in \mathbf{GL}(2n,k) \ \bigg| \  
\begin{matrix}
	\text{$Z_{ij}\in \mathbf{Mat}(2,k)$} \\ \text{is a diagonal matrix}\\
	\text{for every $i,j\in \{1,\dots, n\}$}
\end{matrix}
\right\}.
\]
We define a map
\begin{align}\label{A:GLnDsplit}
	f_{D} : \mathbf{GL}(n,D) & \longrightarrow \mathbf{GL}(n,k) \times \mathbf{GL}(n,k) \\
	Z &\longmapsto (Z^{(1)}, Z^{(2)})
\end{align}
by setting $Z^{(1)}$ to be the matrix obtained from $Z$ by taking the upper-left entry of each $2 \times 2$ block matrix in $Z$,
and $Z^{(2)}$ to be the matrix obtained from $Z$ by taking the lower-right entry of each of the $2 \times 2$ blocks in $Z$.
It is easy to check that $f_{D}$ is an isomorphism of algebraic $k$-groups.

In both cases, $\mathbf{GL}(n,D)$ is a $k$-subgroup of 
$\mathbf{GL}(2n,k)$ of dimension $2n^2$, 
which will appear as a symmetric subgroup inside 
$\mathbf{GL}(n,C)$ in the situations considered later.

\subsection{Special linear groups.}

The group $\mathbf{SL}(n,\HH)$ is the group of $\R$-rational points of the inner form of type $A_{2n-1}$. More generally, we have the following structural result.

\begin{lem}\label{L:splitsemisimple}
	Let $k$ be a field of characteristic $\neq 2$ and let $C$ be an associative composition algebra over $k$.
	We set $G:=\mathbf{SL}(n,C)$.
	Then the following assertions hold: 
	\begin{enumerate}
		\item $G$ is a connected, semisimple and simply-connected algebraic $k$-group.
		\item If $C$ is split (i.e. $C\cong\mathbf{Mat}(2,k)$ when $\dim_k C=4$), then $G\cong\mathbf{SL}(2n,k)$ is \emph{split} of Dynkin type $A_{2n-1}$.
		\item If $C$ is a quaternion division algebra, then $G$ is the \emph{inner $k$-form} of type $A_{2n-1}$ which splits over the quadratic splitting field of $C$. Over $\R$ its group of real points is the Lie group $\mathbf{SU}^{\!*}(2n)$ (compact for $n=1$, non-compact of real rank $n-1$ for $n\ge 2$).
		\item If $\dim_k C=2$, then we have 
		\[
		\begin{aligned}
			C &= k \times k
			&&\Longrightarrow\quad
			G \;\cong\; \mathbf{SL}(n,k)\times\mathbf{SL}(n,k), \\[6pt]
			C &\text{ a quadratic field } L/k
			&&\Longrightarrow\quad
			G \;\cong\; \operatorname{Res}_{L/k}\mathbf{SL}(n,L),
		\end{aligned}
		\]
		so that $G$ has Dynkin type $A_{n-1}\times A_{n-1}$ in the split commutative case and inner type~$A_{n-1}$ in the field case.
		(Here, $\operatorname{Res}_{L/k}$ stands for the restriction of scalars functor.)
	\end{enumerate}
\end{lem}
\begin{proof}
	The split case follows from the identification $C\cong\mathbf{Mat}(2,k)$ and the Cayley-Dickson embedding $\mathbf{GL}(n,C)\cong\mathbf{GL}(2n,k)$, whence $G$ coincides with the commutator subgroup $\mathbf{SL}(2n,k)$.
	
	Suppose $C$ is a quaternion division algebra.  The main theorem of Dieudonn\'e~\cite{Dieudonne1943} shows that $\mathbf{SL}(n,C)$ is the derived subgroup of $\mathbf{GL}(n,C)$, hence semisimple.  Let $L/k$ be the quadratic splitting field of~$C$.  Base-change yields $G_L\cong\mathbf{SL}(2n,L)$, which is split and simply connected, so $G$ is an inner $k$-form of type $A_{2n-1}$ and is itself simply connected.
	
	The two-dimensional cases are direct: for $C=k\times k$, the determinant-one condition gives $\mathbf{S}(\mathbf{GL}(n,k)\times\mathbf{GL}(n,k))\cong\mathbf{SL}(n,k)\times\mathbf{SL}(n,k)$, whereas when $C=L$ is a quadratic field extension, $G\cong\operatorname{Res}_{L/k}\mathbf{SL}(n,L)$.
	
	Connectedness of $G$ is immediate as it is the derived subgroup of a connected reductive group.  This completes the proof.
\end{proof}

The group $\mathbf{SL}(n,\HH)$  is sometimes denoted $\mathbf{SU}^{\!*}(2n)$ in the theory of real
Lie groups.  We record for later use that
\[
\operatorname{rank}_{\R}\mathbf{SL}(n,\HH)=n-1,\qquad
\dim_{\R}\mathbf{SL}(n,\HH)=4n^{2}-1.
\]

Let $G$ be a linear algebraic group defined over $k$, and let
$\bar{k}$ denote the algebraic closure of $k$.
We say that $G$ is reductive if the base change
$G_{\bar{k}} := G \otimes_k \bar{k}$ contains no nontrivial unipotent normal connected linear algebraic subgroup defined over $k$. 
The most typical example of such a group is the general linear group $\mathbf{GL}(n,k)$ (for some $n\geq 1$) defined over $k$.
We are interested in the homogeneous spaces of general linear groups with entries not from $k$ but from an associative composition algebra $C$ over $k$.
On one hand, by using the symplectic representation of $G:=\mathbf{GL}(n,C)$, where $C$ is an associative composition algebra over $k$, it is easy to see that $G$ is a reductive $k$-group. 
On the other hand, all homogeneous spaces of $G$ that we are concerned with can be expressed as quotient of a smaller group.  
A linear algebraic $k$-group is called semisimple if the base change
$G_{\bar{k}} := G \otimes_k \bar{k}$  contains no nontrivial solvable connected linear algebraic subgroup defined over $k$. 
It is well-known that the derived (commutator) subgroup of a reductive $k$-group is always a semisimple $k$-group (\cite[Proposition 14.2 (2)]{Borel}).
In particular, the derived subgroup of a general linear group $\mathbf{GL}(n,k)$ (except when $n=2$ and $k$ has 2 or 3 elements) is the special linear group 
$\mathbf{SL}(n,k)$ consisting of $n\times n$ matrices with entries from $k$ and determinant 1. 
For noncommutative composition algebras $C$, 
the definition of a determinant on $\mathbf{GL}(n,C)$ 
requires additional care. 
In this context one uses the reduced norm, 
which coincides with the usual determinant 
after base change to a splitting field. 
An accessible discussion of determinants over quaternion algebras 
can be found in~\cite{Aslaksen1996}.
\medskip

We continue with the understanding that $G$ denotes $\mathbf{GL}(n,C)$, where $C$ is a quaternion algebra over $k$. 
\begin{defn}
	For $n\geq 2$, the $n$-th special linear group with entries from $C$, denoted $\mathbf{SL}(n,C)$, is the derived subgroup of $\mathbf{GL}(n,C)$. 
\end{defn}
According to Dieudonn\'e~\cite[Theorem 1]{Dieudonne1943}, if $C$ is a skew-field, the commutator subgroup of $\mathbf{GL}(n,C)$ satisfies 
$$
\mathbf{GL}(n,C) / \mathbf{SL}(n,C) \cong C^*/ [C^*,C^*].
$$
Hence, the center of $\mathbf{GL}(n,C)$ is given by the center of the multiplicative group of the skew-field $C$. 

\medskip

\begin{ex}
	In this example, we consider the Hamilton quaternion algebra as the underlying skew-field, $C=\HH$. 
	Let us identify $\mathbf{Mat}(n,\HH)$ with the image of its symplectic representation $\sigma$.
	Then the {\em Study determinant} of an element $Z:=A+jB \in \mathbf{Mat}(n,\HH)$ is defined as follows: 
	$$
	S.\det (Z) := | \det (\sigma(Z)) |^2,
	$$
	where $ \det (\sigma(Z))$ is the ordinary determinant of the complex $2n\times 2n$ matrix $\begin{bmatrix} A & - \overline{B} \\B & A \end{bmatrix}$.
	The Study determinant $S.\det (Z)$ is a surjective group homomorphism 
	$$
	S.\det : \mathbf{GL}(n,\HH) \to \R_+^\times,
	$$
	where $\R_+^\times$ is the multiplicative group of positive real numbers.
	Since $\R_+^\times$ is commutative,  one has
\[
[\mathbf{GL}(n,\HH),\mathbf{GL}(n,\HH)] \subseteq \ker(S.\det).
\]
Moreover, $\ker(S.\det)$ coincides with $\mathbf{SL}(n,\HH)$,
and the group $\mathbf{SL}(n,\HH)$ is perfect.
It follows that
\[
[\mathbf{GL}(n,\HH),\mathbf{GL}(n,\HH)] = \ker(S.\det).
\]
In particular,
\begin{equation}\label{A:SLn}
\mathbf{SL}(n,\HH)=\ker(S.\det).
\end{equation}
	Since the Study determinant is a polynomial map in real coordinates, $\mathbf{SL}(n,\HH)$ is a one-codimensional real algebraic subvariety of $\mathbf{Mat}(n,\HH)$:
	$$
	\dim ( \mathbf{SL}(n,\HH)) = 4n^2 -1. 
	$$
	The theorem of Dieudonn\'e shows that (\ref{A:SLn}) can be used for defining special linear group over $\HH$. 
\end{ex}

\subsection{Symmetric spaces of associative composition algebras.}

Recall that a symmetric space of $G$ is a homogeneous space of the form $G/H$, where $G^\sigma \subseteq H \subseteq N_G(G^\sigma)$
for some automorphism $\sigma : G\to G$ such that $\sigma^2$ is the trivial automorphism of $G$. 
We will prove Proposition~\ref{intro:P1} from the introductory section. 
Let us recall its statement:
\medskip

Let $C$ be an associative composition algebra over a field $k$. 
Let $D$ be a composition subalgebra of $C$ such that $C$ is obtained from $D$ by the Cayley-Dickson doubling process.
For $n\geq 1$, let $G$ (resp. $H$) denote $\mathbf{GL}(n,C)$ (resp. $\mathbf{GL}(n,D)$). Then, after base change to a splitting field $S$ of $C$,
the homogeneous space
\[
G_S/H_S
\]
is a symmetric space.

\medskip

\begin{proof}

    Assume that $C$ is non-split.
Let $L\subset C$ be a maximal subfield; equivalently, $L$ is a quadratic field
extension of $k$.
Fix a splitting field $S/k$ of $C$ containing $L$.
Then
\[
C_S:=C\otimes_k S \cong \mathbf{Mat}(2,S),
\qquad
D_S:=D\otimes_K S \cong S\times S,
\]
where $D_S$ identifies with the diagonal subalgebra of $\mathbf{Mat}(2,S)$.
(Concretely, under the standard matrix realization coming from the Cayley--Dickson
presentation, the subfield $L=k(\sqrt a)$ embeds diagonally as in
\eqref{A:CDdoubling12}, and after base change $L\hookrightarrow K$ this yields
the diagonal identification of $D_S$.)

Consequently,
\[
G_S=\mathbf{GL}(n,C_S)\cong \mathbf{GL}(n,\mathbf{Mat}(2,S))
\cong \mathbf{GL}(2n,S),
\]
and similarly
\[
H_S=\mathbf{GL}(n,D_S)\cong \mathbf{GL}(n,S\times S)
\cong \mathbf{GL}(n,S)\times \mathbf{GL}(n,S),
\]
embedded block-diagonally in $\mathbf{GL}(2n,S)$.

Let
\[
J=\begin{pmatrix} I_n & 0 \\ 0 & -I_n\end{pmatrix}\in \mathbf{GL}(2n,S),
\]
and define an involutive automorphism
\[
\sigma:\mathbf{GL}(2n,S)\longrightarrow \mathbf{GL}(2n,S),
\qquad
\sigma(g)=JgJ^{-1}.
\]
Then $\sigma^2=\mathrm{id}$ and
\[
\mathbf{GL}(2n,S)^\sigma
=
\left\{
\begin{pmatrix} A & 0 \\ 0 & B\end{pmatrix}
\ \middle|\ A,B\in \mathbf{GL}(n,S)
\right\}
\cong \mathbf{GL}(n,S)\times \mathbf{GL}(n,S).
\]
Under the identifications above, this fixed point subgroup coincides with $H_S$.
Therefore $G_S/H_S$ is a symmetric space.
	
	Next, we assume that $C$ is split, that is, $\mathbf{Mat}(2,k)$. Thus we may take the splitting field $S$ to be $k$ itself.
	In this case also, $D$ can be identified with the diagonal algebra. 
	The corresponding groups, as we showed in the previous section, are given by $\mathbf{GL}(2n,k)$ and $\mathbf{GL}(n,k)\times \mathbf{GL}(n,k)$. It is well-known that this subgroup is a symmetric Levi subgroup of $\mathbf{GL}(2n,k)$. 
	Hence, the proof of our assertion is complete. 
\end{proof}

\subsection{Quaternionic symplectic group.}

We will use the symplectic representation of $\mathbf{GL}(n,\HH)$ to define a subgroup which will be instrumental in our cohomology computations. 
The {\em quaternionic symplectic group}, denoted $\mathbf{Sp}(n)$, is defined as the intersection of the image of the symplectic representation of $\mathbf{GL}(n,\HH)$ 
and the $2n$-th unitary group, that is, 
$$
\mathbf{Sp}(n) : = \mathbf{U}(2n,\C) \cap \sigma( \mathbf{GL}(n,\HH)).
$$
When it is clear from the context, we will view $\mathbf{Sp}(n)$ as a subgroup of $\mathbf{GL}(n,\HH)$ by identifying it with the subgroup 
$\sigma^{-1}(\mathbf{U}(2n,\C))$.

The group $\mathbf{Sp}(n)$ is a compact, connected Lie group 
and coincides with the compact real form of 
$\mathbf{SL}(n,\mathbb{H})$.

\section{A Quaternionic Homogeneous Space}\label{S:Quaternionic}

Let ${\rm Gr}(k,\HH^n)$ denote the Grassmann manifold of all $k$-dimensional quaternionic subspaces of $\HH^n$.
Every ordered basis of $\HH^n$ is mapped onto another ordered basis by an invertible matrix $F\in \mathbf{GL}(n,\HH)$, and this matrix $F$ can be made a quaternionic-unitary matrix. 
This means that ${\rm Gr}(k,\HH^n)$ is a homogeneous space for both $\mathbf{GL}(n,\HH)$ and its subgroup $\mathbf{Sp}(n) \subseteq \mathbf{GL}(n,\HH)$. 
On one hand, we have the stabilizer group of a point $W\in {\rm Gr}(k,\HH^n)$ in $\mathbf{Sp}(n)$, which is isomorphic to $\mathbf{Sp}(k)\times \mathbf{Sp}(n-k)$, where $\mathbf{Sp}(k)$ acts on $W$ and $\mathbf{Sp}(n-k)$ acts on the orthogonal complement $W^\perp$ with respect to the quaternionic Hermitian inner product. 
On the other hand, the stabilizer of $W$ in $\mathbf{GL}(n,\HH)$ is a `maximal parabolic subgroup' $\mathbf{P}(k,\HH)$ consisting of matrices of the form 
\begin{align*}
	\begin{bmatrix}
		A & B \\
		0 & C
	\end{bmatrix} \in \mathbf{GL}(n,\HH), \quad{\text{ where }} A \in \mathbf{GL}(k,\HH)\ \text{ and }\ C \in \mathbf{GL}(n-k,\HH).
\end{align*}
We now have two quotients that are isomorphic to each other as real algebraic varieties: 
$$
\mathbf{GL}(n,\HH) / \mathbf{P}(k,\HH) \  \ \cong \  \  \mathbf{Sp}(n) / (\mathbf{Sp}(k)\times \mathbf{Sp}(n-k)).
$$

Let $\tau$ be an involutive antiholomorphic automorphism of an irreducible complex affine variety $M$.
If the set $M_0:= M^\tau$ of fixed points of $\tau$ contains at least one simple (nonsingular) point, then $M_0$ is said to be a {\em real form} of $M$. 
In the case of $M:=\mathbf{GL}(2n,\C)$, the map 
$$
Z\mapsto - E_{2n} \overline{Z} E_{2n}, \quad\text{where}\quad E_{2n}:=\begin{bmatrix} 0 & I_n \\ -I_n & 0 \end{bmatrix},
$$
gives an involutive antiholomorphic automorphim of $\mathbf{GL}(2n,\C)$. 
The fixed point subgroup is precisely the image of the general linear group of quaternions $\mathbf{GL}(n,\HH)$ under its symplectic representation (\ref{A:definingrep}).

\medskip

Let $G$ be a complex linear algebraic group. A closed subgroup $H\subseteq G$ is said to be {\em parabolic} if the quotient variety $G/H$ is a complex projective algebraic variety. 
In this case, the Lie subalgebra of $H$ in ${\rm Lie}(G)$ is called a {\em parabolic Lie subalgebra}. 
Now, let us assume that $G$ is a real linear algebraic group, and $H$ is a real algebraic subgroup of $G$.  
Let us denote the Lie algebra of $H$ (resp. of $G$) by $\mathfrak{h}$ (resp. by $\mathfrak{g}$). 
In this notation, $H$ is called a {\em parabolic subgroup} of $G$ if and only if the complexification of $\mathfrak{h}$ 
is a parabolic subalgebra of the complexification of the Lie algebra $\mathfrak{g}$.

\medskip

Let $n\geq 2$. 
Evidently, the complex general linear group $\mathbf{GL}(n,\C)$ is a real algebraic subgroup of $\mathbf{GL}(n,\HH)$.
To see this, recall that every element $Z\in \mathbf{GL}(n,\HH)$ has the form $Z=A+ j B$, where $A$ and $B$ are from $\mathbf{Mat}(n,\C)$. 
By requiring $B=0$, we obtain precisely the matrices in $\mathbf{GL}(n,\C)$, showing that it is a closed (in Zariski topology) subgroup of $\mathbf{GL}(n,\HH)$. 

We continue to work with the image of the symplectic representation $\sigma$ of $\mathbf{GL}(n,\HH)$. 
Then $\mathbf{GL}(n,\C)$ is identified with the subgroup consisting of block diagonal matrices, 
$$
\sigma(\mathbf{GL}(n,\C)) = 
\left\{ \begin{bmatrix} A &0 \\ 0 & \overline{A} \end{bmatrix} \mid A \in \mathbf{GL}(n,\C) \right\}.
$$
Its Lie algebra is given by 
$$
{\rm Lie}(\sigma(\mathbf{GL}(n,\C))) = 
\left\{ \begin{bmatrix} A &0 \\ 0 & \overline{A} \end{bmatrix} \mid A \in \mathbf{Mat}(n,\C) \right\}.
$$
Let us denote this Lie subalgebra by $W$ to simplify our notation. 
Clearly, $\overline{W} = W$.
Hence, $W$ is already a complex Lie subalgebra of ${\rm Lie}(\mathbf{GL}(2n,\C))$.  
However, $W$ can not be a parabolic subalgebra for dimension reasons.
Indeed, a conjugate of a parabolic subalgebra must contain the Lie algebra of all upper triangular matrices in ${\rm Lie}(\mathbf{GL}(2n,\C))$.
But the dimension of the algebra of upper triangular matrices in ${\rm Lie}(\mathbf{GL}(2n,\C))$ is bigger than the dimension of $W$.
We summarize these observations as follows.

\begin{lem}\label{L:notaparabolic}
	The real algebraic subgroup $H:=\mathbf{GL}(n,\C)$ in $G:=\mathbf{GL}(n,\HH)$ is {\em not} a parabolic subgroup. 
\end{lem}

Although Lemma~\ref{L:notaparabolic} states that $\mathbf{GL}(n,\C)$ is not a parabolic subgroup of $\mathbf{GL}(n,\HH)$,
we observe that it comes very close to being one. 
Let $  \mathbf{P}_{n,2n}$ denote the (complex) parabolic subgroup of $\mathbf{GL}(2n,\C)$, defined by 
$$
\mathbf{P}_{n,2n} = \left\{ \begin{bmatrix} A &B \\ 0 & C \end{bmatrix} \in \mathbf{GL}(2n,\C) \mid A,C \in \mathbf{GL}(n,\C) \right\}.
$$
In particular, $\mathbf{GL}(2n,\C)/  \mathbf{P}_{n,2n}$ is the Grassmann variety ${\rm Gr}(n,\C^{2n})$. 
The image of $\mathbf{GL}(n,\C)$ in $\mathbf{GL}(2n,\C)$ under the symplectic representation is contained in the parabolic subgroup $  \mathbf{P}_{n,2n}$.
In fact, it is easy to check that the preimage $\sigma^{-1}( \mathbf{P}_{n,2n})$ is exactly the subgroup $\mathbf{GL}(n,\C)$ in $\mathbf{GL}(n,\HH)$. 
This means that the image of the compositions 
\begin{align}\label{A:projectGLnintoGr}
	\mathbf{GL}(n,\HH)\xrightarrow{\sigma} \mathbf{GL}(2n,\C) \xrightarrow{\pi} {\rm Gr}(n,\C^{2n}),
\end{align}
where $\pi$ is the canonical quotient morphism, is precisely the homogeneous space 
$$
\mc{Y}_n:=\mathbf{GL}(n,\HH) / \mathbf{GL}(n,\C).
$$
Equivalently, this means that the $\mathbf{GL}(n,\HH)$-orbit of the origin in ${\rm Gr}(n,\C^{2n})$, that is, 
the $\mathbf{GL}(n,\HH)$-orbit of $1 \mathbf{P}_{n,2n} / \mathbf{P}_{n,2n}$, is isomorphic to $\mc{Y}_n$. 
We count dimensions: 
\begin{align*}
	\dim_\C {\rm Gr} ( n, \C^{2n}) &= n (2n-n) = n^2 \\
	\dim_\R \mc{Y}_n &= \dim_\R \mathbf{GL}(n,\HH) -\dim_\R \mathbf{GL}(n,\C) = 4n^2 - 2n^2 = 2n^2.
\end{align*}
Hence, the dimensions of ${\rm Gr} ( n, \C^{2n})$, considered as a real algebraic variety, and the real algebraic homogeneous space $\mc{Y}_n$ are equal.

\begin{prop}
	The complex Grassmann variety ${\rm Gr}(n,\C^{2n})$ is a $\mathbf{GL}(n,\HH)$-equivariant, smooth completion of the real algebraic homogeneous space $\mc{Y}_n$.
\end{prop}
\begin{proof}
	We already observed that the homogeneous space $\mc{Y}_n$ is isomorphic to the $\mathbf{GL}(n,\HH)$-orbit of the point 
	$1 \mathbf{P}_{n,2n} / \mathbf{P}_{n,2n}$ in ${\rm Gr}(n,\C^{2n})$. 
	By Lemma~\ref{L:notaparabolic}, $\mathbf{GL}(n,\C)$ is not a parabolic subgroup of $\mathbf{GL}(n,\HH)$, implying that $\mc{Y}_n$ is not a complete variety. 
	Additionally, we observed earlier that $\dim_\R \mc{Y}_n = \dim_\R {\rm Gr}(n,\C^{2n})$. 
	But a maximal dimensional orbit of a (real) algebraic group action is always an open subset of the variety. 
	Hence, we conclude that $\mc{Y}_n$ is an open real algebraic subvariety of the Grassmannian, which is a smooth real algebraic variety,
	and $\mc{Y}_n\neq {\rm Gr}(n,\C^{2n})$.
	This finishes the proof. 
\end{proof}

\begin{rem}
	It would be very interesting to determine all orbit types of $\mathbf{GL}(n,\HH)$ in ${\rm Gr}(n,\C^{2n})$. 
\end{rem}

\subsection{Equivariant Cohomology.}

In this subsection we will discuss the (equivariant) cohomologies of our homogeneous spaces.

We begin with $\mc{Y}_n := \mathbf{GL}(n,\HH)/\mathbf{GL}(n,\C)$. 
The following observation that will simplify much of our work.

The following well-known result~\cite[Proposition 1.143]{Knapp} will be useful:
	\medskip
	
	\begin{lem}\label{L:fromKnapp}
		Let $G\subseteq \mathbf{GL}(m,\C)$ be a closed linear group that is the common zero locus of some set of real-valued polynomials in the real and imaginary parts of the matrix entries, and let $\mathfrak{g}$ be its linear Lie algebra. 
		Suppose that $G$ is closed under taking adjoints. Let $K$ be the group $G\cap \mathbf{U}(m)$, and let $\mathfrak{p}$ be the subspace of Hermitian matrices in $\mathfrak{g}$.
		Then the map $K\times \mathfrak{p} \to G$ given by $(k,X) \mapsto k e^X$ is a diffeomorphism. 
	\end{lem}
	To give it a name, we call the decomposition of $G$ in the previous lemma the `polar decomposition' of $G$. 
	In Knapp's book, it is stated that the map $(k,X) \mapsto k e^X$ is a homeomorphism, in this setting the exponential maps on vector groups are diffeomorphisms. 
\begin{prop}\label{P:Yn}
Let $G=\sigma(\GL(n,\HH))$ and $H=\sigma(\GL(n,\C))$.
Let $K=\Sp(n)$ and $M=U(1)\cdot SU(n)$ be maximal compact subgroups
of $G$ and $H$, respectively.
Then there is a $K$-equivariant diffeomorphism
\[
G/H \;\cong\; K \times_M (\mathfrak p/\mathfrak p'),
\]
where $\mathfrak p/\mathfrak p'$ is the $M$-representation
coming from the Cartan decomposition.
In particular, $G/H$ is $K$-equivariantly homotopy equivalent to
$K/M$, and hence
\[
\GL(n,\HH)/\GL(n,\C)
\simeq
\Sp(n)/(U(1)\cdot SU(n)).
\]
\end{prop}

\begin{proof}
	
	The general linear group $\mathbf{GL}(n,\HH)$ has a positive dimensional (real) center. 
	Recall that the center of $\mathbf{GL}(n,\HH)$ is contained in $\mathbf{GL}(n,\C)$. 
	It follows that 
	$$
	\mathbf{GL}(n,\HH)/\mathbf{GL}(n,\C) \cong \mathbf{SL}(n,\HH)/\left(\mathbf{SL}(n,\HH) \cap \mathbf{GL}(n,\C)\right).
	$$
	We look closely at the intersection $\mathbf{SL}(n,\HH) \cap \mathbf{GL}(n,\C)$, which is explicitly described via the symplectic representation: 
	\begin{align*}
		\sigma(\mathbf{SL}(n,\HH) \cap \mathbf{GL}(n,\C)) &= \sigma(\mathbf{SL}(n,\HH)) \cap \sigma(\mathbf{GL}(n,\C))\\
		&=\left\{ \begin{bmatrix} A & 0 \\ 0 & \overline{A} \end{bmatrix} \in \mathbf{GL}(2n,\C) \mid  (\det A)(\det \overline{A}) = 1 \right\} \\
		&=  \left\{ \begin{bmatrix} A & 0 \\ 0 & \overline{A} \end{bmatrix} \in \mathbf{GL}(2n,\C) \mid  \det A  \in \mathbf{U}(1) \right\} \\
		&=  \left\{  \begin{bmatrix} x \cdot A & 0 \\ 0 & \bar{x} \cdot\overline{A} \end{bmatrix} \in \mathbf{GL}(2n,\C) \mid  \det A =1,\  x\in \mathbf{U}(1)\right\}.
	\end{align*}
	Hence, we have the following identification of homogeneous spaces, 
	\begin{align}\label{A:identification}
		\mathbf{GL}(n,\HH) / \mathbf{GL}(n,\C) \ \ \cong \ \ \mathbf{SL}(n,\HH) / \mathbf{SLU}(n,\C),
	\end{align}
	where $\mathbf{SLU}(n,\C)$ is the subgroup 
	$$
	\mathbf{SLU}(n,\C):=  \mathbf{U}(1)\mathbf{SL}(n,\C),
	$$
	which is an internal direct product of $\mathbf{U}(1)$ with $\mathbf{SL}(n,\C)$.
	
	\medskip
	
	We apply the polar decomposition to the symplectic representation of $\mathbf{SL}(n,\HH)$ and its subgroup $\mathbf{SLU}(n,\C)$.
	By intersecting $\sigma(\mathbf{SL}(n,\HH))$ with $\mathbf{U}(2n)$ we obtain exactly $\mathbf{Sp}(n)$. 
	By intersecting $\sigma(\mathbf{SLU}(n,\C))$ with $\mathbf{U}(2n)$ we obtain: 
	$$
	\mathbf{U}(2n) \cap \sigma(\mathbf{SLU}(n,\C)) \cong \mathbf{U}(1)\times \mathbf{SU}(n).
	$$
	Then we get the diagram of smooth maps in Figure~\ref{F:polar}, where $\Phi$ (and its restriction $\Phi'$) is a diffeomorphism, and the vertical maps are natural embeddings.
	\begin{figure}[htp]
		\centering
		\begin{tikzpicture}[>=stealth, scale=1.5]
			\node (Spn) at (0,0) {$\mathbf{Sp}(n) \times \mathfrak{p}$};
			\node (GLnHH) at (3,0) {$\sigma(\mathbf{SL}(n,\mathbb{H}))$};
			\node (Un) at (0,-1) {$\mathbf{U}(1)\times \mathbf{SU}(n) \times \mathfrak{p}'$};
			\node (GLnC) at (3,-1) {$\sigma(\mathbf{SLU}(n,\C))$};
			\draw[->] (Spn) -- node[above] {$\Phi$} (GLnHH);
			\draw[<-] (Spn) -- node[left] { } (Un);
			\draw[->] (Un) -- node[below] {$\Phi'$} (GLnC);
			\draw[<-] (GLnHH) -- node[right] { } (GLnC);
		\end{tikzpicture}
		\caption{Polar decompositions.}
		\label{F:polar}
	\end{figure}
   
It follows from the polar decompositions  applied to
$G:=\sigma(\mathbf{SL}(n,\HH))$ and $H:=\sigma(\mathbf{SLU}(n,\C))$ that the quotient
$G/H$ is $K$-equivariantly diffeomorphic to the homogeneous vector bundle
over $K/M$ associated with the $M$-representation $\mathfrak p/\mathfrak p'$,
namely
\[
G/H \;\cong\; K \times_M (\mathfrak p/\mathfrak p'),
\]
where $K=G\cap U(2n)=\mathbf{Sp}(n)$ and $M=H\cap U(2n)=U(1)\cdot SU(n)$.
Combining this with the identification~\eqref{A:identification} yields
\[
\mathbf{GL}(n,\HH)/\mathbf{GL}(n,\C)\ \cong\ \mathbf{Sp}(n)\times_{\,U(1)\cdot SU(n)}(\mathfrak p/\mathfrak p').
\]
Finally, the fiberwise contraction $[k,v]\mapsto [k,tv]$ defines a $K$-equivariant
strong deformation retraction of $\mathbf{Sp}(n)\times_{M}(\mathfrak p/\mathfrak p')$ onto the
zero section $\mathbf{Sp}(n)/M$. Hence $\mathbf{GL}(n,\HH)/\mathbf{GL}(n,\C)$ is homotopy equivalent to
$\mathbf{Sp}(n)/(U(1)\cdot SU(n))$, completing the proof.
\end{proof}

\begin{rem}
	According to Berger~\cite[Tableau II, Case 17 (in pg. 158)]{Berger1957}, the homogeneous space
	$\mathbf{SL}(n,\HH) / \mathbf{SLU}(n,\C)$ is a symmetric space. 
	Note: Berger uses $\mathbf{SU}^*(2n)$ to denote our $\mathbf{SL}(n,\HH)$.
	Additionally, he uses $\mathbf{T}$ to denote the unitary group $\mathbf{U}(1)$. 
\end{rem}

The following corollary follows from our previous proposition.

\begin{cor}
There is an isomorphism of graded $\Q$-algebras
\[
H^*\!\left(\mathbf{GL}(n,\HH)/\mathbf{GL}(n,\C);\Q\right)
\cong
H^*\!\left(\mathbf{Sp}(n)/( \mathbf{U}(1)\cdot \mathbf{SU}(n));\Q\right).
\]
\end{cor}

\begin{proof}
By Proposition~\ref{P:Yn}, the homogeneous space
$\mathbf{GL}(n,\HH)/\mathbf{GL}(n,\C)$
is diffeomorphic to the associated vector bundle
\[
\mathbf{Sp}(n)\times_{\,\mathbf{U}(1)\cdot \mathbf{SU}(n)}
(\mathfrak p/\mathfrak p').
\]
The fiberwise contraction gives a deformation retraction onto the
zero section $\mathbf{Sp}(n)/(\mathbf{U}(1)\cdot \mathbf{SU}(n))$.
The result follows.
\end{proof}

We now recall a classical theorem recorded in~\cite[Theorem 26.1]{Borel1953}.
Let $U$ be a closed subgroup of a compact Lie group $G$ having the same rank as $G$. 
This means that a maximal torus of $G$, denoted $T$, is entirely contained in $U$.
Let 
$$
{\rm rank}\ G  = {\rm rank}\ U = \dim T = l.
$$
Let $W(G)$ and $W(U)$ denote the corresponding Weyl groups of $(G,T)$ and $(U,T)$, respectively.
Then the Poincar\'e polynomial of $G/U$ is given by 
\begin{align}\label{A:Hirschformula}
	P(G/U,t) = \frac{ (1-t^{2s_1}) (1-t^{2s_2})\cdots (1-t^{2s_l})}  { (1-t^{2r_1}) (1-t^{2r_2})\cdots (1-t^{2r_l})} .  
\end{align}
where $s_1,\dots, s_l$, resp. $r_1,\dots, r_l$, are the degrees of the fundamental $W(G)$-invariant polynomial functions on ${\rm Lie}(T)$, resp. $W(U)$-invariant polynomials, on the algebra of polynomial functions on ${\rm Lie}(T)$.
This formula follows from~\cite[Theorem 26.1 and Proposition 27.2 (b)]{Borel1953}.
In our case, that is $(G,U)=(\mathbf{Sp}(n),\mathbf{U}(1)\mathbf{SU}(n))$, the degrees of the fundamental invariants are easy to find.
On one hand, since the Weyl group of $\mathbf{Sp}(n)$ is a Weyl group of type $C_n$, 
we have the degrees $2,4,6,\dots, 2n$. 
On the other hand, we have the following decomposition of the Weyl group of ${\mathbf{U}(1)\mathbf{SU}(n)}$:
$$
W(\mathbf{U}(1)\mathbf{SU}(n)) \cong W(\mathbf{U}(1)) \times W(\mathbf{SU}(n)).
$$
Since $W(\mathbf{SU}(n))$ is a Weyl group of type $A_{n-1}$, its degrees are given by $2,3,\cdots,n$.
The  group $W(\mathbf{U}(1))$ is trivial.
Therefore, the degrees of the fundamental invariants of $W(\mathbf{U}(1)\mathbf{SU}(n))$
are given by $1,2,3,4,\dots, n$.
Then the Poincar\'e polynomial of $\mathbf{Sp}(n)/\mathbf{U}(1)\mathbf{SU}(n)$ (hence, of $\mc{Y}_n$) is given by 
\begin{align}
	P(\mathbf{Sp}(n)/\mathbf{U}(1)\mathbf{SU}(n),t) &= \frac{(1-t^{4}) (1-t^{8})(1-t^{12})\cdots (1-t^{4n})}{(1-t^{2})(1-t^{4}) (1-t^{6})\cdots (1-t^{2n})} \notag \\
	&= \prod_{i=1}^n (1+t^{2i}). \label{A:PofY_n}
\end{align}
For example, if $n=3$, then we have
\begin{align*}
	P(\mathbf{Sp}(3)/\mathbf{U}(1)\mathbf{SU}(3),t) &= \frac{(1-t^{4}) (1-t^{8})(1-t^{12})}{(1-t^{2}) (1-t^{4})(1-t^{6})} \\
	&= (1+t^2)(1+t^4) (1+t^6).
\end{align*}
\medskip

Our next goal is to compute the equivariant cohomology of $\mc{Y}_n$ with respect to the subgroup
$\mathbf{U}(1)^n \subset \mathbf{U}(1) \mathbf{SU}(n) \subset \mathbf{Sp}(n)$.
\medskip

We begin with a general observation.

\begin{lem}\label{L:bundle_retract}
Let $K$ be a compact Lie group, $M\subseteq K$ a closed subgroup,
and $V$ a finite-dimensional real $M$-representation.
Set $E:=K\times_M V$.
Then the zero section $K/M\hookrightarrow E$ is a $K$-equivariant strong deformation retract.
In particular, for any closed subgroup $T\subseteq K$ and any coefficient ring $R$,
the inclusion induces an isomorphism of graded rings
\[
H_T^*(E;R)\cong H_T^*(K/M;R).
\]
\end{lem}
\begin{proof}
Define $H:[0,1]\times E\to E$ by $H(t,[k,v])=[k,tv]$.
This is well-defined because scalar multiplication commutes with the $M$-action on $V$.
It is $K$-equivariant and gives a strong deformation retraction onto the zero section.
Equivariant cohomology is invariant under equivariant homotopy equivalence, hence the claim.
\end{proof}

In Proposition~\ref{P:Yn} we identified $\mc{Y}_n=\mathbf{GL}(n,\HH)/\mathbf{GL}(n,\C)$
(with its $\mathbf{U}(1)^n$-action by left multiplication) with the
homogeneous vector bundle
\[
\mathbf{Sp}(n)\times_{\,\mathbf{U}(1)\cdot\mathbf{SU}(n)}(\mathfrak p/\mathfrak p'),
\]
arising from the polar decomposition diffeomorphism
\begin{equation}\label{A:contractionmap}
\mathbf{Sp}(n)\times \mathfrak{p} \longrightarrow \mathbf{GL}(n,\HH),\qquad (k,X)\longmapsto k e^X,
\end{equation}
and its analogue for $\mathbf{GL}(n,\C)$.
Here $\mathfrak p$ (resp.\ $\mathfrak p'$) denotes the space of Hermitian matrices
in the Lie algebra of $\GL(n,\HH)$ (resp.\ $\mathbf{GL}(n,\C)$).
The action of $\mathbf{U}(1)^n$ on the associated bundle is induced from the
left action on the $\Sp(n)$-factor.
By Lemma~\ref{L:bundle_retract} (fiberwise contraction), the zero section
$\Sp(n)/(\mathbf{U}(1)\cdot\mathbf{SU}(n))$ is a $\mathbf{U}(1)^n$-equivariant
deformation retract of $\mc{Y}_n$. Hence we obtain:

\begin{cor}\label{C:eqcoHH/C}
	The $\mathbf{U}(1)^n$-equivariant cohomology of $\mathbf{GL}(n,\HH)/\mathbf{GL}(n,\C)$ is given by 
	$$
	H^*_{\mathbf{U}(1)^n}( \mathbf{GL}(n,\HH)/\mathbf{GL}(n,\C);\Q) \cong H^*_{\mathbf{U}(1)^n}( \mathbf{Sp}(n)/\mathbf{U}(1) \mathbf{SU}(n);\Q).
	$$
\end{cor}

\subsection{The equivariant cohomology of \texorpdfstring{$\mc{Z}_n$}.}

In this subsection, we consider our second homogeneous space $\mc{Z}_n := \mathbf{GL}(2n,\R)/\mathbf{GL}(n,\C)$.
We notice the isomorphism, $\mc{Z}_n \cong \mathbf{SL}(2n,\R)/\mathbf{SLU}(n,\C)$.
This is the noncompact symmetric space that appears in~\cite[Tableau II, Case 13 (in pg. 158)]{Berger1957}.
We will discuss its homology. 
The proof of the following proposition (and of Proposition~\ref{P:Wn}) follows from general results in~\cite{Mostow1962}. We provide direct proofs for the convenience of readers.

\begin{prop}\label{P:Zn}
	The noncompact symmetric space $\mc{Z}_n$ has a deformation retract onto the compact homogeneous space 
	$\mathbf{SO}(2n,\R) / \mathbf{U}(1)\mathbf{SU}(n)$.
\end{prop}

\begin{proof}
	
	We apply Lemma~\ref{L:fromKnapp} again. 
	By intersecting $\mathbf{SL}(2n,\R)$ with $\mathbf{U}(2n)$ we obtain exactly $\mathbf{SO}(2n,\R)$. 
	By intersecting $\mathbf{SLU}(n,\C)$ with $\mathbf{U}(2n)$ we obtain $\mathbf{U}(1)  \mathbf{SU}(n)$.
	Then we get a diagram of smooth maps as in Figure~\ref{F:polar2}, where $\Phi$ (and its restriction $\Phi'$) is a diffeomorphism, and the vertical maps are natural embeddings.
	\begin{figure}[htp]
		\centering
		\begin{tikzpicture}[>=stealth, scale=1.5]
			\node (Spn) at (0,0) {$\mathbf{SO}(2n,\R) \times \mathfrak{p}$};
			\node (GLnHH) at (3,0) {$\sigma(\mathbf{SL}(2n,\R))$};
			\node (Un) at (0,-1) {$\mathbf{U}(1)  \mathbf{SU}(n) \times \mathfrak{p}'$};
			\node (GLnC) at (3,-1) {$\sigma(\mathbf{SLU}(n,\C))$};
			\draw[->] (Spn) -- node[above] {$\Phi$} (GLnHH);
			\draw[<-] (Spn) -- node[left] { } (Un);
			\draw[->] (Un) -- node[below] {$\Phi'$} (GLnC);
			\draw[<-] (GLnHH) -- node[right] { } (GLnC);
		\end{tikzpicture}
		\caption{Polar decomposition of $\mathbf{GL}(2n,\R)/\mathbf{GL}(n,\C)$.}
		\label{F:polar2}
	\end{figure}
	It follows that 
	\begin{align*}
		\mathbf{SO}(2n,\R) / \mathbf{U}(1) \mathbf{SU}(n) \times (\mathfrak{p} / \mathfrak{p}') \xrightarrow{\resizebox{!}{.8em}{$\sim$}}  \mathbf{GL}(2n,\R) /\mathbf{GL}(n,\C).
	\end{align*}
	 
    By Lemma~\ref{L:fromKnapp}, applied to $G=\mathbf{SL}(2n,\R)$ and
$H=\mathbf{SLU}(n,\C)$, we obtain polar decompositions
\[
G \cong K \exp(\mathfrak p),
\qquad
H \cong M \exp(\mathfrak p'),
\]
where
\[
K=\mathbf{SO}(2n,\R),
\qquad
M=\mathbf{U}(1)\cdot \mathbf{SU}(n).
\]
It follows that the homogeneous space $G/H$ is
diffeomorphic to the associated vector bundle
\[
\mathbf{SO}(2n,\R)
\times_{\,M}
(\mathfrak p/\mathfrak p').
\]
Since scalar multiplication in the fibers defines a
$\mathbf{SO}(2n,\R)$-equivariant deformation retraction
onto the zero section
$\mathbf{SO}(2n,\R)/M$,
the space $\mc{Z}_n$ deformation retracts onto
$\mathbf{SO}(2n,\R)/M$.
\end{proof}

By the associated-bundle description in Proposition~\ref{P:Zn} and the equivariant deformation retraction of homogeneous vector bundles (Lemma~\ref{L:bundle_retract}), we obtain the following.

\begin{cor}\label{C:eqco2nR/nC}
    There is an isomorphism of graded $\Q$-algebras
\[
H^*_{\mathbf{U}(1)^n}\!\left(
\mathbf{GL}(2n,\R)/\mathbf{GL}(n,\C);\Q
\right)
\;\cong\;
H^*_{\mathbf{U}(1)^n}\!\left(
\mathbf{SO}(2n)/(\mathbf{U}(1)\cdot\mathbf{SU}(n));\Q
\right).
\]
\end{cor}

Next, we determine the Betti numbers of our symmetric space.
Since $\mathbf{SO}(2n,\R)$ and $\mathbf{U}(1) \mathbf{SU}(n)$ have equal ranks, we apply, as before, Hirsch's formula, proved in~\cite[Theorem 26.1 and Proposition 27.2 (b)]{Borel1953}.

The Weyl group of $\mathbf{SO}(2n,\R)$ is a Weyl group of type $D_n$.
Its degrees are well-known to be $2,4,6,\dots, 2n-2,n$. 
Earlier, we found the degrees of the fundamental invariants of $W( \mathbf{U}(1)\mathbf{SU}(n) )$. 
They are given by $2,3,4,\dots, n$.
It follows that the Poincar\'e polynomial of $\mathbf{SO}(2n,\R)/\mathbf{U}(1)\mathbf{SU}(n)$ is given by 
\begin{align*}
	P(\mathbf{SO}(2n,\R)/\mathbf{U}(1)\mathbf{SU}(n),t) &= \frac{(1-t^{4}) (1-t^{8})(1-t^{12})\cdots (1-t^{4n-4})(1-t^{2n})}{(1-t^2)(1-t^{4}) (1-t^{6})\cdots (1-t^{2n})}\\
	&= \prod_{i=1}^{n-1} (1+t^{2i}).
\end{align*}
For example, for $n=3$, we have
\begin{align*}
	P(\mathbf{SO}(6,\R)/\mathbf{U}(1)\mathbf{SU}(3),t) &= \frac{(1-t^{4}) (1-t^{8})(1-t^{6})}{(1-t^{2}) (1-t^{4})(1-t^{6})} \\
	&=(1+t^2) (1+t^4).
\end{align*}
\medskip

\subsection{The equivariant cohomology of \texorpdfstring{$\mc{W}_n$}.}

In this subsection, we consider our last homogeneous space $\mc{W}_n = \mathbf{GL}(n,\HH_s)/\mathbf{GL}(n,\C_s)$.
We showed earlier that 
$$
\mathbf{GL}(n,\C_s) \cong \mathbf{GL}(n,\R)\times \mathbf{GL}(n,\R)\quad\text{and}\quad 
\mathbf{GL}(n,\HH_s) \cong \mathbf{GL}(2n,\R). 
$$
Thus, after taking out the center of $\mathbf{GL}(2n,\R)$, we have 
$$
\mc{W}_n \cong \mathbf{SL}(2n,\R)/ \R(\mathbf{SL}(n,\R)\times \mathbf{SL}(n,\R)).
$$
This is the noncompact symmetric space that appears in~\cite[Tableau II, Case 10 (in pg. 158)]{Berger1957}.

\begin{prop}\label{P:Wn}
	The noncompact symmetric space $\mc{W}_n$ has a deformation retract onto the compact homogeneous space 
	$\mathbf{SO}(2n,\R) / \mathbf{SO}(n,\R)\times \mathbf{SO}(n,\R)$.
\end{prop}

\begin{proof}
	We apply Lemma~\ref{L:fromKnapp} once more. 
	By intersecting $\mathbf{SL}(2n,\R)$ with $\mathbf{U}(2n)$, we obtain exactly $\mathbf{SO}(2n,\R)$. 
	By intersecting the subgroup $\R(\mathbf{SL}(n,\R)\times \mathbf{SL}(n,\R))$ with $\mathbf{U}(2n)$, we get the product subgroup $\mathbf{SO}(n,\R)\times \mathbf{SO}(n,\R)$.
Hence, as in the proofs of Propositions~\ref{P:Yn} and~\ref{P:Zn}, we obtain an identification
of homogeneous spaces of the form
\[
\mathbf{GL}(2n,\R)/\mathbf{GL}(n,\C_s)
\;\cong\;
\mathbf{SO}(2n,\R)\times_{\,\mathbf{SO}(n,\R)\times\mathbf{SO}(n,\R)}(\mathfrak p/\mathfrak p'),
\]
where $\mathfrak p/\mathfrak p'$ is the $\mathbf{SO}(n,\R)\times\mathbf{SO}(n,\R)$-module
arising from the corresponding Cartan decompositions.
By fiberwise contraction, the associated vector bundle deformation retracts
onto the zero section $\mathbf{SO}(2n,\R)/(\mathbf{SO}(n,\R)\times\mathbf{SO}(n,\R))$,
which proves the assertion. 
\end{proof}

Next, we apply Lemma~\ref{L:bundle_retract} to the associated bundle description above
(with $K=\mathbf{SO}(2n,\R)$ and $M=\mathbf{SO}(n,\R)\times\mathbf{SO}(n,\R)$).
This yields the following. 

\begin{cor}\label{C:eqco2nR/nCs}
	There is an isomorphism of graded $\Q$-algebras
\[
H^*_{\mathbf{U}(1)^n}\!\left(\mathbf{GL}(2n,\R)/\mathbf{GL}(n,\C_s);\Q\right)
\cong
H^*_{\mathbf{U}(1)^n}\!\left(\mathbf{SO}(2n,\R)/(\mathbf{SO}(n,\R)\times \mathbf{SO}(n,\R));\Q\right).
\]
\end{cor}

The compact homogeneous space $\mathbf{SO}(2n,\R)/(\mathbf{SO}(n,\R)\times \mathbf{SO}(n,\R))$ is the Grassmann manifold of oriented $n$-planes in $\R^{2n}$,
which is another classical object studied by a number of authors. 
The structure of its cohomology space, in particular, its Betti numbers depend on the parity of $n$. 
A detailed description of its cohomology ring is recorded in the textbook~\cite[Table III, pg. 495-496]{GHV}. 
The Poincar\'e polynomial of this Grassmann manifold is given by the $t^2$-analog of the binomial numbers, $\begin{bmatrix} n \\ \lfloor n/ 2 \rfloor \end{bmatrix}_{t^2}$, where $ \lfloor n/ 2 \rfloor$ stands for the largest integer less than $n/2$. 

\subsection{Proof of Theorem~\ref{intro:1}.}

We recall the statement of Theorem~\ref{intro:1}:
\medskip

Let $D$ and $C$ be composition algebras in $\{\C,\C_s,\HH,\HH_s\}$
such that there exists an injective algebra homomorphism $D\hookrightarrow C$. Let $T:=\mathbf{U}(1)^n$.

Then the following assertions hold.

\begin{enumerate}

\item[(1)] If $(C,D)=(\HH,\C)$, then
$
P(\mc{Y}_n,t)=\prod_{i=1}^{n}(1+t^{2i}).
$
Moreover, there is an isomorphism of graded $S$-algebras
\[
H_T^*(\mc{Y}_n;\Q)
\;\cong\;
S \otimes_{S^{W(\mathbf{Sp}(n))}}
S^{W(\mathbf{U}(1)\cdot \mathbf{SU}(n))}.
\]

\item[(2)] If $(C,D)=(\HH_s,\C)$, then
$
P(\mc{Z}_n,t)=\prod_{i=1}^{n-1}(1+t^{2i}).
$
Moreover, there is an isomorphism of graded $S$-algebras
\[
H_T^*(\mc{Z}_n;\Q)
\;\cong\;
S \otimes_{S^{W(\mathbf{SO}(2n))}}
S^{W(\mathbf{U}(1)\cdot \mathbf{SU}(n))}.
\]

\item[(3)] If $(C,D)=(\HH_s,\C_s)$, then
$
P(\mc{W}_n,t)=
\begin{bmatrix}
n \\ \lfloor n/2\rfloor
\end{bmatrix}_{t^2}
$ (Gaussian binomial coefficient).
Moreover, there is an isomorphism of graded $S$-algebras
\[
H_T^*(\mc{W}_n;\Q)
\;\cong\;
S \otimes_{S^{W(\mathbf{SO}(2n))}}
S^{W(\mathbf{SO}(n)\times \mathbf{SO}(n))}.
\]

\end{enumerate}

\begin{proof}[Proof of Theorem~\ref{intro:1}]
We begin with case \textup{(1)}.
The Poincar\'e polynomial of $\mc{Y}_n$ was computed in~\eqref{A:PofY_n}.
By Corollary~\ref{C:eqcoHH/C}, it remains to determine the
$T$-equivariant cohomology of
\[
\Sp(n)/(\mathbf{U}(1)\cdot \mathbf{SU}(n)),
\qquad T:=\mathbf{U}(1)^n.
\]

Let $G:=\Sp(n)$ and $H:=\mathbf{U}(1)\cdot \mathbf{SU}(n)$.
We work with rational coefficients and set
\[
S:=H_T^*(pt;\Q)\cong \Q[x_1,\dots,x_n], \qquad \deg x_i=2.
\]

By~\cite[\S1, Proposition~1(iii)]{Brion1997}, there is an isomorphism of graded
$S$-algebras
\begin{equation}\label{A:combine1}
H_T^*(G/H;\Q)
\;\cong\;
S \otimes_{S^{W(G)}} H_G^*(G/H;\Q),
\end{equation}
where $W(G)$ denotes the Weyl group of $(G,T)$.

Next, by~\cite[\S1, Remarks~3), Examples~1), and Proposition~1(i)]{Brion1997},
\begin{equation}\label{A:combine2}
H_G^*(G/H;\Q)
\;\cong\;
H_H^*(pt;\Q)
\;\cong\;
S^{W(H)},
\end{equation}
where $W(H)$ is the Weyl group of $(H,T)$.

Since $W(G)\cong W(\Sp(n))$ is the Weyl group of type $C_n$
(the hyperoctahedral group),
and $W(H)\cong W(\mathbf{U}(1)\cdot \mathbf{SU}(n))\cong S_n$,
combining~\eqref{A:combine1} and~\eqref{A:combine2}
yields
\[
H_T^*(\mc{Y}_n;\Q)
\;\cong\;
S \otimes_{S^{W(G)}} S^{W(H)},
\]
which proves the second assertion of \textup{(1)}.

The proofs of \textup{(2)} and \textup{(3)} are analogous.
Using Corollaries~\ref{C:eqco2nR/nC} and~\ref{C:eqco2nR/nCs},
together with the corresponding compact symmetric pairs
\[
(\mathbf{SO}(2n),\,\mathbf{U}(1)\cdot\mathbf{SU}(n))
\quad\text{and}\quad
(\mathbf{SO}(2n),\,\mathbf{SO}(n)\times\mathbf{SO}(n)),
\]
and applying the same argument as above,
we obtain the claimed descriptions of their $T$-equivariant cohomology rings.
\end{proof}

\section{Homogeneous Spaces Related to Clifford Algebras}\label{S:Final}

We wish to highlight an interesting and nontrivial case of the opening question of our paper.  In the remainder of this section, the base field will be taken as either $F=\R$  or $F=\C$.
A useful contemporary reference for the basics of this section (under the assumption that $F\in \{\R,\C\}$) is the book by Vaz and Rold\~ao, \cite{vaz2016}.
\medskip

Let $M$ be an $F$-vector space and of dimension $n$ and $Q$ be a quadratic form on $M$. We will call $(M,Q)$ an $F$-quadratic space.
We assume that the characteristic of $F$ is not 2. 
Let $T(M)$ be the full tensor algebra over $M$.
Let $I(M)$ denote the ideal generated in $T(M)$ by the elements $x\otimes x - Q(x)\cdot 1$, where $x\in M$. 
The factor algebra $T(M)/I(M)$ is called the {\rm Clifford algebra} of the pair $(M,Q)$ and denoted ${\rm Cl}(M,Q)$. 
In particular, $M$ is a subspace of ${\rm Cl}(M,Q)$. 
Denote the unit group of ${\rm Cl}(M,Q)$ by $G({\rm Cl}(M,Q))$.
It acts on ${\rm Cl}(M,Q)$ by conjugation. 
Then the following subgroup of the stabilizer of $M$ in $G({\rm Cl}(M,Q))$ is called the {\em Clifford group} of the pair $(M,Q)$:
$$
\Gamma(M,Q) = \{ g\in G({\rm Cl}(M,Q)) \mid g M g^{-1} \subseteq M,\ Q(gmg^{-1}) = Q(m) \text{ for all $m\in M$}\}.
$$
The subgroup of $\Gamma(M,Q)$ consisting of elements of the form $v_{1}\cdots v_{2k}$, where $Q(v_i) = \pm 1$ for $1\leq i \leq 2k$, is called the {\em spin group}, denoted ${\rm Spin}(M,Q)$. 
Equivalently, we have
$$
{\rm Spin}(M,Q)= {\rm Pin}(M,Q) \ \cap \  {\rm Cl}^+(M,Q),
$$
where 
$$
{\rm Pin}(M,Q) := \{ g\in \Gamma(M,Q) \mid g Q(g) = 1\}
$$
and 
${\rm Cl}^+(M,Q) \subset {\rm Cl}(M,Q)$ is the subalgebra  spanned by the elements that can be written as a product of even number of generators.
It is worth mentioning that there are exact natural sequences:
\[1\rightarrow F^\times \rightarrow \Gamma(M,Q) \rightarrow {\mathbf O}(M)\rightarrow 1\]
and
\[1\rightarrow \mathbb{Z}/2\Z \rightarrow {\rm Spin}(M,Q) \rightarrow \mathbf{SO}(M)\rightarrow 1\]
where $\mathbf{O}(M)$ (resp. $\mathbf{SO}(M)$) is the orthogonal group (resp. special orthogonal group) on $M$ preserving the bilinear form associated with $Q$.

It is well-known since Chevalley's work~\cite{Chevalley} on the Clifford algebras that ${\rm Cl}(M,Q)$ is either a central simple algebra or a direct sum of two central simple algebras depending on the dimension $\dim M$. 
Hence, by the Weddernburn-Artin theorem, ${\rm Cl}(M,Q)$ is isomorphic to a full matrix algebra over a quaternion algebra $C$, or to a direct sum of two copies of such an algebra. 
Accordingly, we view the Clifford and the spin groups as subgroups of the unit group of the corresponding matrix algebra over $C$. 
Let $D\subseteq C$ be a maximal associative composition subalgebra of $C$.
We then get a homogeneous space $\Gamma(M,Q)/\Gamma(M,Q)_D$ (resp. ${\rm Spin}(M,Q)/{\rm Spin}(M,Q)_D$), where 
$\Gamma(M,Q)_D$ (resp. ${\rm Spin}(M,Q)_D$) is the subgroup consisting of $D$-rational
points of $\Gamma(M,Q)$ (resp. of ${\rm Spin}(M,Q)$).

In the case $(M,Q)$ is a real quadratic space, if we denote by $(M_{\mathbb{C}}, Q_{\mathbb{C}})$ its complexification, then
\[{\rm Cl}(M_{\mathbb{C}}, Q_{\mathbb{C}})\cong{\rm Cl}_{\mathbb{C}}(M,Q)\]
Therefore, the complex Clifford Algebra structure follows from the study of the real case. Fix $(M,Q)$ a real quadratic space where $M$ has dimension $n$ and $Q$ is a non-degenerate quadratic form of signature $(p,q)$ ($p+q=n$). We denote the Clifford Algebra of $(M,Q)$ by ${\rm Cl}(p,q)$. It is known the characterization of these ones in terms of matrix algebras, as mentioned before, and these are determined by $r=p-q\, \mod\,8$: (where \(m:=\lfloor n/2\rfloor\))
\[
\mathrm{Cl}(p,q)\cong
\begin{cases}
\mathbf{Mat}(2^m,\R), & r\equiv 0,2 \pmod 8,\\
\mathbf{Mat}(2^m,\R)\oplus\mathbf{Mat}(2^m,\R), & r\equiv 1 \pmod 8,\\
\mathbf{Mat}(2^m,\C), & r\equiv 3,7 \pmod 8,\\
\mathbf{Mat}(2^{m-1},\HH), & r\equiv 4,6 \pmod 8,\\
\mathbf{Mat}(2^{m-1},\HH)\oplus\mathbf{Mat}(2^{m-1},\HH), & r\equiv 5 \pmod 8.
\end{cases}
\]

\begin{ex}
	${\rm Cl}(0,1)\cong \mathbb{C}$, ${\rm Cl}(1,0)\cong \mathbb{R}\oplus \mathbb{R}=\mathbb{C}_s$, ${\rm Cl}(0,2)\cong \mathbb{H}$, ${\rm Cl}(1,1)\cong{\rm Cl}(2,0)\cong \mathbf{Mat}(2,\mathbb{R})=\mathbb{H}_s$
\end{ex}

The complex Clifford Algebra depends only on the parity of $n = p + q$, denote by ${\rm Cl}(n,\mathbb{C})$ the complexification of ${\rm Cl}(p,q)$, thus:

\begin{table}[h!]
	\centering
	\begin{tabular}{|c|c|}
		\hline
		\(n=2k\) even &  \({\rm Cl}(2k,\mathbb{C})\cong \mathbf{Mat}(2^{k}, \mathbb{C}) \)         \\ \hline
		\(n=2k+1\) odd &  \({\rm Cl}(2k+1,\mathbb{C})\cong \mathbf{Mat}(2^{k}, \mathbb{C})\oplus \mathbf{Mat}(2^{k}, \mathbb{C})  \) \\ \hline
	\end{tabular}
\end{table}

\medskip

Similarly, we will denote by $\Gamma(p,q)$ the real Clifford group and by $\Gamma(n,\mathbb{C})$ the complex Clifford group, such that $\Gamma(n,\mathbb{C})$ is the complexification of $\Gamma(p,q)$. For $F=\mathbb{R}$, ${\rm Spin}(p,q)$ is a double cover of the identity component of $\mathbf{SO}(p,q)$ and if $p\geq 2$ or $q\geq 2$, it is connected.
For $F=\mathbb{C}$, ${\rm Spin}(n)$ is a double cover of $\mathbf{SO}(n,\mathbb{C})$, it is connected and simply connected, that is, it is a universal cover of $\mathbf{SO}(n,\mathbb{C})$. The real group ${\rm Spin}(n,0)$ is the maximal compact subgroup of ${\rm Spin}(p,q)$.

Due to the adjoint or vector representation $\sigma: \Gamma(p,q)\to \rm{Aut}{\rm Cl}(p,q) $ defined by $\sigma(A)(X)=AXA^{-1}$, which image is an orthogonal group with $\mathrm{Ker}(\sigma)=\mathbb{R}^{\ast}$, we have the following isomorphism:
$$\Gamma(p,q)/\mathbb{R}^\ast\cong \left \{ \begin{array}{cc}
	\mathbf{O}(p,q),   & p+q \quad even \\
	\mathbf{SO}(p,q),   & p+q\quad odd
\end{array}\right.
$$
Taking the complexification we get 
$$\Gamma(n,\mathbb{C})/\mathbb{C}^\ast\cong \left \{ \begin{array}{cc}
	\mathbf{O}(n,\mathbb{C}),   & n \quad even \\
	\mathbf{SO}(n, \mathbb{C}),   & n\quad odd
\end{array}\right.
$$
Considering these isomorphisms, we can now study involutions of $\Gamma(n)$ whose fixed points are the real Clifford group $\Gamma(p,q)$. Involutions of the orthogonal group $\mathbf{O}(n,\mathbb{C})$ or $\mathbf{SO}(n,\mathbb{C})$, depending on the parity of $n$, give rise to involutions of $\Gamma(n,\mathbb{C})$, as $F^\ast$ is the identity through the adjoint representation.

Let 
\[
I_{p,q} = \begin{pmatrix} I_p & 0 \\ 0 & -I_q \end{pmatrix}.
\]
An example of such involution is the automorphism 
$\theta: \mathbf{SO}(n,\mathbb{C}) \rightarrow \mathbf{SO}(n,\mathbb{C})$ defined by $\theta(X)=I_{p,q} \bar{X}^{-t}I_{p,q}$. The fixed points are given by $\mathbf{SO}(p,q)$.

Now, we would like to discuss the (equivariant) cohomology of the homogeneous space $\Gamma(n,\mathbb{C})/\Gamma(p,q)$. In order to do that, we start by studying the homogeneous spaces $\mc{O}_n:=\mathbf{O}(n,\mathbb{C})/\mathbf{O}(p,q)$ and
$\mc{SO}_n:=\mathbf{SO}(n,\mathbb{C})/\mathbf{SO}(p,q)$, this last one is the noncompact symmetric space that appears in~\cite[Tableau II, Case 31 (in pg. 158)]{Berger1957}.

Let $K$ denote either $\mathbf{O}(n,\R)$ or $\mathbf{SO}(n,\R)$,
depending on the space under consideration.
Let $T\subset \mathbf{SO}(n,\R)\subset K$
be a maximal torus.
Then $T\cong \mathbf{U}(1)^r$ with $r=\lfloor n/2\rfloor$.
Via the natural inclusion $K\hookrightarrow \mathbf{SO}(n,\C)$ we regard $T$ as acting on $\mc{O}_n$ and $\mc{SO}_n$ by left multiplication.
In the sequel, all equivariant cohomology groups will be taken with respect to this torus $T$.

\begin{prop}\label{P:SO}
	The noncompact symmetric spaces $\mc{O}_n$ and $\mc{SO}_n$ have a $K$-equivariant deformation retract onto the compact homogeneous spaces 
	$\mathbf{O}(n,\R) /  (\mathbf{O}(p)\times \mathbf{O}(q))$ and  $\mathbf{SO}(n,\R) /  (\mathbf{SO}(p)\times \mathbf{SO}(q))$ respectively, where $K=\mathbf{O}(n,\R)$ in the first case and
$K=\mathbf{SO}(n,\R)$ in the second.
\end{prop}

\begin{proof} 
We prove the statement for $\mc{SO}_n=\mathbf{SO}(n,\C)/\mathbf{SO}(p,q)$;
the argument for $\mc{O}_n$ is identical.
View $\mathbf{SO}(n,\C)$ as a real Lie group. Its maximal compact subgroup is
$K:=\mathbf{SO}(n,\R)$.
Applying Lemma~\ref{L:fromKnapp} to $G=\mathbf{SO}(n,\C)$, we obtain a polar
decomposition
\[
\Phi:K\times \mathfrak p \longrightarrow \mathbf{SO}(n,\C),\qquad (k,X)\mapsto k\,e^{X},
\]
where $\mathfrak p$ is the subspace of Hermitian elements in the Lie algebra of $G$.
Similarly, applying Lemma~\ref{L:fromKnapp} to the real form
$H:=\mathbf{SO}(p,q)\subset \mathbf{SO}(n,\C)$, we obtain
\[
\Phi':L\times \mathfrak p' \longrightarrow \mathbf{SO}(p,q),
\qquad (l,Y)\mapsto l\,e^{Y},
\]
where $L:=H\cap K=\mathbf{SO}(p)\times \mathbf{SO}(q)$ and $\mathfrak p'$ is the
corresponding Hermitian subspace for $H$.

Passing to quotients yields a $K$-equivariant diffeomorphism
\[
K\times_{L}(\mathfrak p/\mathfrak p')
\;\xrightarrow{\ \sim\ }\;
\mathbf{SO}(n,\C)/\mathbf{SO}(p,q)=\mc{SO}_n,
\]
where $L$ acts on $\mathfrak p/\mathfrak p'$ via the isotropy representation.
Since $\mathfrak p/\mathfrak p'$ is a real vector space, it is contractible, and
the fiberwise contraction
\[
[k,v]\longmapsto [k,tv],\qquad t\in[0,1],
\]
defines a $K$-equivariant deformation retraction of
$K\times_{L}(\mathfrak p/\mathfrak p')$ onto the zero section, which identifies
with $K/L=\mathbf{SO}(n,\R)/(\mathbf{SO}(p)\times\mathbf{SO}(q))$.
This proves that $\mc{SO}_n$ admits a $K$-equivariant deformation retraction onto
$\mathbf{SO}(n,\R)/(\mathbf{SO}(p)\times\mathbf{SO}(q))$.
The case of $\mc{O}_n$ is proved in the same way, replacing special orthogonal
groups by orthogonal groups throughout.
\end{proof}

By Proposition~\ref{P:SO}, the spaces $\mc{O}_n$ and $\mc{SO}_n$
admit $T$-equivariant deformation retractions onto the corresponding
compact homogeneous spaces.
Since equivariant cohomology is invariant under equivariant homotopy equivalence,
we obtain the following.

\begin{cor}\label{C:eqco2nR/nClast}
The following isomorphisms of graded $\Q$-algebras hold:
\begin{enumerate}
\item
$$
H_T^*( \mc{O}_n;\Q)\ \cong\
H_T^*\!\left(\mathbf{O}(n,\R)/(\mathbf{O}(p)\times\mathbf{O}(q));\Q\right).
$$
\item
$$
H_T^*( \mc{SO}_n;\Q)\ \cong\
H_T^*\!\left(\mathbf{SO}(n,\R)/(\mathbf{SO}(p)\times\mathbf{SO}(q));\Q\right).
$$
\end{enumerate}
\end{cor}

Next, we determine the Betti numbers of our symmetric spaces.
It is well known that
\[
\mathbf{O}(n,\R)/(\mathbf{O}(p)\times \mathbf{O}(q))
\;\cong\;
\mathrm{Gr}(p,\R^{n}),
\qquad (q=n-p),
\]
the real Grassmannian of $p$-planes in $\R^{n}$.
Hence its Poincar\'e polynomial is given by the $t^2$–analog of the binomial
coefficient:
\begin{align}\label{PPSO}
P\!\left(\mathbf{O}(n,\R)/(\mathbf{O}(p)\times \mathbf{O}(q)),t\right)
&=
\binom{n}{p}_{t^2}
=
\prod_{i=1}^{p}\frac{1-t^{2(q+i)}}{1-t^{2i}}.
\end{align}

Now consider the compact homogeneous space in the odd-dimensional case,
namely $n=2m+1$ with $p=2k$ and $q=n-p=(2m+1)-2k$:
\[
\mathbf{SO}(2m+1,\R)
/\bigl(\mathbf{SO}(2k)\times \mathbf{SO}(2m+1-2k)\bigr).
\]
By~\cite[Table II (p.~494)]{GHV}, its Poincar\'e polynomial is
\begin{align}\label{PPO}
P\!\left(
\mathbf{SO}(2m+1,\R)/
(\mathbf{SO}(2k)\times \mathbf{SO}(2m+1-2k)),
t
\right)
&=
\frac{\displaystyle\prod_{i=m-k+1}^{m} (1 - t^{4i})}
{\displaystyle\left(\prod_{i=1}^{k-1} (1 - t^{4i})\right)(1 - t^{2k})}.
\end{align}

    Thus, we can now characterize the $\mathbf{U}(1)^r$-equivariant cohomology
and the Poincar\'e polynomial of $\Gamma(n,\mathbb{C})/\Gamma(p,q)$.

The scalar subgroups satisfy $\C^\ast/\R^\ast\simeq S^1$ (by polar decomposition
$\C^\ast\simeq S^1\times \R_{>0}$ and $\R^\ast\simeq \{\pm 1\}\times \R_{>0}$).

\begin{thm}\label{clifford:1}
Let $r=\lfloor n/2\rfloor$ and set $q=n-p$.
Then the following statements hold.

\begin{enumerate}

\item The $\mathbf{U}(1)^r$-equivariant cohomology of
$\Gamma(n,\mathbb{C})/\Gamma(p,q)$ with $\Q$-coefficients is given by

\begin{enumerate}
\item[(a)] If $n$ is even,
\[
H^*_{\mathbf{U}(1)^r}\!\left(
\Gamma(n,\mathbb{C})/\Gamma(p,q);\Q
\right)
\cong
H^*(S^1;\Q)\otimes
H^*_{\mathbf{U}(1)^r}\!\left(
\mathbf{O}(n,\R)/(\mathbf{O}(p)\times\mathbf{O}(q));\Q
\right).
\]

\item[(b)] If $n$ is odd,
\[
H^*_{\mathbf{U}(1)^r}\!\left(
\Gamma(n,\mathbb{C})/\Gamma(p,q);\Q
\right)
\cong
H^*(S^1;\Q)\otimes
H^*_{\mathbf{U}(1)^r}\!\left(
\mathbf{SO}(n,\R)/(\mathbf{SO}(p)\times\mathbf{SO}(q));\Q
\right).
\]
\end{enumerate}

\item The Poincar\'e polynomial of $\Gamma(n,\mathbb{C})/\Gamma(p,q)$ is

\begin{enumerate}
\item[(a)] If $n$ is even,
\[
P(\Gamma(n,\mathbb{C})/\Gamma(p,q),t)
=
(1+t)\,
\prod_{i=1}^{p}
\frac{1-t^{2(q+i)}}{1-t^{2i}}.
\]

\item[(b)] If $n=2m+1$ is odd with $p=2k$ and
$q=(2m+1)-2k$,
\[
P(\Gamma(n,\mathbb{C})/\Gamma(p,q),t)
=
(1+t)\,
\frac{\displaystyle\prod_{i=m-k+1}^{m}(1-t^{4i})}
{\displaystyle\left(\prod_{i=1}^{k-1}(1-t^{4i})\right)(1-t^{2k})}.
\]
\end{enumerate}

\end{enumerate}
\end{thm}

\begin{proof}
We treat the case where $n$ is even; the odd case is entirely analogous,
replacing $\mathbf{O}$ by $\mathbf{SO}$ throughout.

\begin{enumerate}

\item
Since $n$ is even, we have
\[
\Gamma(n,\C)/\C^\ast \cong \mathbf{O}(n,\C).
\]
Moreover, the real subgroup $\Gamma(p,q)$ maps onto $\mathbf{O}(p,q)$
with kernel $\R^\ast$.
It follows that
\[
\Gamma(n,\C)/\Gamma(p,q)
\;\cong\;
(\C^\ast/\R^\ast) \times \mc{O}_n,
\qquad
\mc{O}_n=\mathbf{O}(n,\C)/\mathbf{O}(p,q).
\]

Since $\C^\ast/\R^\ast\simeq S^1$, we obtain
\[
\Gamma(n,\C)/\Gamma(p,q)
\;\simeq\;
S^1\times \mc{O}_n.
\] By
Proposition~\ref{P:SO}, $\mc{O}_n$ admits a $T$-equivariant deformation
retract onto
\[
\mathbf{O}(n,\R)/(\mathbf{O}(p)\times\mathbf{O}(q)),
\]
we obtain a $T$-equivariant homotopy equivalence
\[
\Gamma(n,\C)/\Gamma(p,q)
\;\simeq_T\;
S^1 \times
\mathbf{O}(n,\R)/(\mathbf{O}(p)\times\mathbf{O}(q)).
\]

Equivariant cohomology is invariant under equivariant homotopy equivalence,
and the Künneth formula (over $\Q$) yields
\[
H_T^*(\Gamma(n,\C)/\Gamma(p,q);\Q)
\cong
H^*(S^1;\Q)
\otimes
H_T^*\!\left(
\mathbf{O}(n,\R)/(\mathbf{O}(p)\times\mathbf{O}(q));
\Q
\right),
\]
which proves part (1).

\item
Since Poincar\'e polynomials are multiplicative under products,
and $P(S^1,t)=1+t$, we obtain
\[
P(\Gamma(n,\C)/\Gamma(p,q),t)
=
(1+t)\,
P\!\left(
\mathbf{O}(n,\R)/(\mathbf{O}(p)\times\mathbf{O}(q)),
t
\right).
\]
Substituting the expression from~\eqref{PPSO} gives the claimed formula.
\end{enumerate}
\end{proof}


\subsection*{Acknowledgment}
Ana Casimiro's work was funded by national funds through the FCT - Fundação para a Ci\^{e}ncia e a Tecnologia, I.P., under the scope of the projects  
UID/00297/2025 (https://doi.org/10.54499/UID/00297/2025) and UID/PRR/00297/2025 \newline
(https://doi.org/10.54499/UID/PRR/00297/2025) (Center for Mathematics and Applications – NOVA Math).
A part of this paper was written while Ferruh \"{O}zbudak was visiting Department of Mathematics, Tulane University, New Orleans. 
He would like to thank the department for the support and hospitality. 
Ferruh \"{O}zbudak gratefully acknowledges a partial support from the Sabanc{\i} University Research Grants, grant no. B.A.CG-24-02873.

\section*{Author Contributions}
All authors contributed equally to the development of the results and the writing of the manuscript. All authors reviewed and approved the final version.

\section*{Data Availability}
No datasets were generated or analyzed during the current study.

\section*{Declarations}

\textbf{Conflict of interest.}
The authors declare that they have no conflict of interest.
\bibliographystyle{plain}
\bibliography{myreferences}

@book {Chevalley,
    AUTHOR = {Chevalley, Claude},
     TITLE = {The algebraic theory of spinors and {C}lifford algebras},
      NOTE = {Collected works. Vol. 2,
              Edited and with a foreword by Pierre Cartier and Catherine
              Chevalley,
              With a postface by J.-P. Bourguignon},
 PUBLISHER = {Springer-Verlag},
ADDRESS = {Berlin},
      YEAR = {1997},
     PAGES = {xiv+214},
      ISBN = {3-540-57063-2},
   MRCLASS = {01A75 (15A66 81R25)},
  MRNUMBER = {1636473},
}

@incollection {Brion1997,
    AUTHOR = {Brion, Michel},
     TITLE = {Equivariant cohomology and equivariant intersection theory},
 BOOKTITLE = {Representation theories and algebraic geometry ({M}ontreal,
              {PQ}, 1997)},
    SERIES = {NATO Adv. Sci. Inst. Ser. C: Math. Phys. Sci.},
    VOLUME = {514},
     PAGES = {1--37},
      NOTE = {Notes by Alvaro Rittatore},
 PUBLISHER = {Kluwer Acad. Publ., Dordrecht},
ADDRESS = {Montreal},
      YEAR = {1998},
      ISBN = {0-7923-5193-2},
   MRCLASS = {14C17 (14L30 14M15 55N91)},
  MRNUMBER = {1649623},
MRREVIEWER = {Dan Edidin},
}

@book{VeldkampSpringer,
  author    = "Springer, T. A. and Veldkamp, F. D.",
  title     = "Octonions, Jordan Algebras and Exceptional Groups",
  series    = "Springer Monographs in Mathematics",
  pages     = "208",
  publisher = "Springer",
address   = "Berlin",
  year      = "2000",
  doi       = "10.1007/978-3-662-12622-6",
  url       = "https://doi.org/10.1007/978-3-662-12622-6"
}

@article{Hurwitz1922,
  author    = "Hurwitz, A.",
  title     = "{\\\"U}ber die Komposition der quadratischen Formen",
  journal   = "Math. Ann.",
  volume    = "88",
  number    = "1-2",
  pages     = "1--25",
  year      = "1922",
  doi       = "10.1007/BF01448439"
}

@article{Jacobson1958,
  author    = "Jacobson, N.",
  title     = "Composition algebras and their automorphisms",
  journal   = "Rend. Circ. Mat. Palermo (2)",
  volume    = "7",
  pages     = "55--80",
  year      = "1958",
  doi       = "10.1007/BF02854388"
}

@misc{KConrad,
  author    = "Conrad, K.",
  title     = "Quaternion Algebras",
  year      = "2011",
  note      = "Lecture notes available at the author's website",
  url       = "https://kconrad.math.uconn.edu/blurbs/ringtheory/quaternionalg.pdf"
}

@book{Borel,
  author    = "Borel, A.",
  title     = "Linear Algebraic Groups",
  edition   = "2",
  series    = "Graduate Texts in Mathematics",
  volume    = "126",
  pages     = "288",
  publisher = "Springer",
  year      = "1991",
  doi       = "10.1007/978-1-4612-0941-6",
  url       = "https://doi.org/10.1007/978-1-4612-0941-6"
}

@article{Aslaksen1996,
  author    = "Aslaksen, H.",
  title     = "Quaternionic determinants",
  journal   = "Math. Intelligencer",
  volume    = "18",
  number    = "3",
  pages     = "57--65",
  year      = "1996",
  doi       = "10.1007/BF03024312"
}

@article{Dieudonne1943,
  author    = "Dieudonné, J.",
  title     = "Les déterminants sur un corps non commutatif",
  journal   = "Bull. Soc. Math. France",
  volume    = "71",
  pages     = "27--45",
  year      = "1943"
}

@book{Knapp,
  author    = "Knapp, A. W.",
  title     = "Lie Groups Beyond an Introduction",
  edition   = "2",
  series    = "Progress in Mathematics",
  volume    = "140",
  pages     = "812",
  publisher = "Birkhäuser Boston, Inc.",
  address   = "Boston, MA",
  year      = "2002"
}

@article{Berger1957,
  author    = "Berger, M.",
  title     = "Les espaces symétriques noncompacts",
  journal   = "Ann. Sci. École Norm. Sup. (3)",
  volume    = "74",
  pages     = "85--177",
  year      = "1957"
}

@article{Borel1953,
  author    = "Borel, A.",
  title     = "Sur la cohomologie des espaces fibrés principaux et des espaces homogénes de groupes de Lie compacts",
  journal   = "Ann. of Math. (2)",
  volume    = "57",
  pages     = "115--207",
  year      = "1953",
  doi       = "10.2307/1969728"
}

@article{Mostow1962,
  author    = "Mostow, G. D.",
  title     = "Covariant fiberings of Klein spaces. II",
  journal   = "Amer. J. Math.",
  volume    = "84",
  pages     = "466--474",
  year      = "1962",
  doi       = "10.2307/2372983"
}

@book{GHV,
  author    = "Greub, W. and Halperin, S. and Vanstone, R.",
  title     = "Connections, Curvature, and Cohomology",
  series    = "Pure and Applied Mathematics",
  volume    = "47-III",
  pages     = "593",
  publisher = "Academic Press",
  address   = "New York-London",
  year      = "1976",
  note      = "Volume III: Cohomology of principal bundles and homogeneous spaces"
}

@book{vaz2016,
  author    = "Vaz Jr, J. and Rocha Jr, R.",
  title     = "An Introduction to Clifford Algebras and Spinors",
  publisher = "Oxford University Press",
address   = {Oxford},
  year      = "2016"
}

\end{document}